\documentclass[smallextended]{svjour3_CorrectedPostsubmission_Part2}

\usepackage{amsfonts}
\usepackage{graphicx}
\usepackage{epstopdf}
\usepackage{amsmath, amssymb, bm}
\usepackage[fixamsmath,disallowspaces]{mathtools}

\RequirePackage{fix-cm}
\smartqed  
\input{tcilatex}
\journalname{Multibody System Dynamics}
\begin{document}

\vspace{-10ex}\title{Screw and Lie Group Theory in Multibody Dynamics}
\subtitle{Recursive Algorithms and Equations of Motion of Tree-Topology Systems}
\author{Andreas M\"uller}

\date{\vspace{-15ex}
\,
}

\institute{A. M\"uller \at Institute of Robotics, Johannes Kepler University, Altenberger Str. 69, 4040 Linz, Austria, a.mueller@jku.at
}
\titlerunning{Screw Theory in MBS Dynamics}
\maketitle

\begin{abstract}

Screw and Lie group theory allows for user-friendly modeling of multibody
systems (MBS) while at the same they give rise to computationally efficient
recursive algorithms. The inherent frame invariance of such formulations
allows for use of arbitrary reference frames within the kinematics modeling
(rather than obeying modeling conventions such as the Denavit-Hartenberg
convention) and to avoid introduction of joint frames. The computational
efficiency is owed to a representation of twists, accelerations, and
wrenches that minimizes the computational effort. This can be directly
carried over to dynamics formulations. In this paper recursive $O\left(
n\right) $ Newton-Euler algorithms are derived for the four most frequently
used representations of twists, and their specific features are discussed.
These formulations are related to the corresponding algorithms that were
presented in the literature. The MBS motion equations are derived in closed
form using the Lie group formulation. One are the so-called 'Euler-Jourdain'
or 'projection' equations, of which Kane's equations are a special case, and
the other are the Lagrange equations. The recursive kinematics formulations
are readily extended to higher orders in order to compute derivatives of the
motions equations. To this end, recursive formulations for the acceleration
and jerk are derived. It is briefly discussed how this can be employed for
derivation of the linearized motion equations and their time derivatives.
The geometric modeling allows for direct application of Lie group
integration methods, which is briefly discussed.

\end{abstract}


\section{Introduction}

The core task in computational multibody system (MBS) dynamics is to either
construct the equations of motion (EOM) explicitly, that can be written for
an unconstrained tree-topology MBS in the form%
\begin{equation}
\mathbf{M}\left( \mathbf{q}\right) \ddot{\mathbf{q}}+\mathbf{C}\left( \dot{%
\mathbf{q}},\mathbf{q}\right) \dot{\mathbf{q}}=\mathbf{Q}\left( \dot{\mathbf{%
q}},\mathbf{q},t\right) ,  \label{EOM}
\end{equation}%
in a way that is easy to pursue, or to evaluate them for given $\left( \ddot{%
\mathbf{q}},\dot{\mathbf{q}},\mathbf{q}\right) $ and $t$, respectively to
solve them, in a computationally efficient way for $\mathbf{q}\left(
t\right) $. In continuation of \cite{Part1} the aim of this paper is to
present established $O\left( n\right) $ formulations in a common geometric
setting and to show that this setting allows for a flexible and
user-friendly MBS modeling.

Screw and Lie group theory provides a geometric framework that allows for
achieving optimal computational performance and at the same time allows for
an intuitive and flexible modeling. In particular, it gives rise to a
formulation of the MBS kinematics that does not involve body-fixed joint
frames. The kinematics modeling is indeed reflected in the formulation used
to evaluate the EOM. A central concept is the representation of velocities
(twists) as screws. Four different variants were recalled in \cite{Part1}.
In this paper their application to dynamics modeling is reviewed. A
well-known approach, which exploits the fact that rigid body twists are
screws, is the so-called 'spatial vector' formulation introduced in \cite%
{Featherstone1983,Featherstone2008}, respectively the so-called 'spatial
operator algebra' that was formalized in \cite{Rodriguez1991}. The latter is
the basis for the $O\left( n\right) $ forward dynamics algorithms introduced
in \cite%
{Fijany1995,Jian1991,Jian1995,LillyOrin1991,Rodriguez1987,Rodriguez1992}.
The fundamental operation underlying these formulations is the frame
transformations of screws, i.e. twists and wrenches. The fact that the
latter can be expressed in terms of compact matrix operations gave rise to a
matrix formulation for the MBS kinematic and dynamics \cite%
{Angeles2003,Legnani1,Legnani2,Uicker2013} using screw algebra. While these
formulations make merely use of the algebraic properties of screws (e.g.
velocities, accelerations, wrenches) several algorithms for generating the
EOM of MBS with tree topology were reported that also exploit the fact that
finite rigid body motions constitute the Lie group $SE\left( 3\right) $
whose Lie algebra $se\left( 3\right) $ is isomorphic to the algebra of
screws \cite{Brockett1984,Hervet1978,Hervet1982,Chevalier1991,Chevalier1994}%
. The central relation is the \emph{product of exponentials} (POE)
introduced in \cite{Brockett1984}. The important feature of such a geometric
Lie group formulation is the frame invariance, which makes it independent
from any modeling convention like Denavit-Hartenberg. This allows for direct
processing of CAD data, and gives further rise to numerically advantageous
Lie group time integration methods. Yet there is no established Lie group
algorithm for the generation respectively evaluation of the EOM that takes
full advantage of the freedom to chose different motion representations
enabled by the frame invariance.

This paper is organized as follows. Recursive relations for the acceleration
and jerk, and thus for the time derivatives of the Jacobians, are first
derived in section \ref{secAcc}. The Newton-Euler equations for the four
different representations of twists introduced in \cite{Part1} are then
recalled in section \ref{secNE}. The corresponding recursive $O\left(
n\right) $ inverse dynamics algorithm for evaluating the EOM are presented
in section \ref{secRecEOM}. The body-fixed algorithm is similar to that in 
\cite%
{AndersonCritchley2003,Bae2001,Fijany1995,Hollerbach1980,HsuAnderson2002,LillyOrin1991,Liu1988,Park1994,ParkBobrowPloen1995,PloenPark1997,ParkKim2000,SohlBobrow2001}%
, the hybrid formulation to that in \cite%
{Anderson1992,BaeHwangHaug1988,Jian1991,Jian1995,Rodriguez1991,Rodriguez1992}%
, and the spatial formulation to that in \cite{Featherstone2008}. Two
versions of the EOM in closed form are presented in section \ref%
{secClosedEOM}. In section \ref{secProjEqu} the 'Euler-Jourdain'
respectively 'projection' equations \cite{BremerBook2010,Wittenburg1977} are
presented that, together with the screw formulation of MBS kinematics,
allows for an efficient MBS modeling in terms of readily available geometric
data. In section \ref{secLagrange} a closed form of the Lagrangian EOM is
presented using the Lie group approach. It should be noticed that the
presented formulations allow for modeling MBS without introduction of joint
frames, while applying the recursive kinematics and dynamics algorithm that
is deemed best suited. The significance of the Lie group formulation for the
linearization of the EOM as well as the determination of derivative of the
EOM w.r.t. geometric design parameters and time derivatives is discussed in
section \ref{secLinearization}. Finally in section \ref{secIntegration} the
application of Lie group integration methods is briefly discussed. The
kinematic relations that were presented in \cite{Part1} are summarized in
appendix A. The basic Lie group background can be found in \cite%
{LynchPark2017,Murray,Selig}.

\section{Acceleration, Jerk, and Partial Derivatives of Jacobian%
\label{secAcc}%
}

Besides the compact description of finite and instantaneous motions of a
system of articulated bodies, a prominent feature of the screw theoretical
approach is that it allows for expressing the partial derivatives explicitly
in terms geometric objects. Moreover, the analytic formulation of the
kinematics using the POE gives rise to compact expressions for higher
derivatives of the instantaneous joint screws, i.e. of the Jacobian, which
may be relevant for sensitivity analysis and linearization of motion
equations. In this section results for the acceleration and jerk of a
kinematic chain are presented for the body-fixed, spatial, and hybrid
representation. The corresponding relations for the mixed representation are
readily found from either one of these using the relations in table 3 of 
\cite{Part1}.

\subsection{Body-Fixed Representation%
\label{secDerBody}%
}

Starting from (\ref{VbJac}) the body-fixed acceleration is $\dot{\mathbf{V}}%
_{i}^{\text{b}}=\mathsf{J}_{i}^{\text{b}}\ddot{\mathbf{q}}+\dot{\mathsf{J}}%
_{i}^{\text{b}}\dot{\mathbf{q}}$, and explicitly in terms of the body-fixed
instantaneous screw coordinates%
\begin{equation}
\dot{\mathbf{V}}_{i}^{\text{b}}=\sum_{j\leq i}\mathbf{J}_{i,j}^{\text{b}}%
\ddot{q}_{j}+\sum_{j\leq i}\sum_{k\leq i}\frac{\partial }{\partial {q_{k}}}%
\mathbf{J}_{i,j}^{\text{b}}\dot{q}_{j}\dot{q}_{k}.  \tag{2}
\end{equation}%
Using the matrix form of (\ref{JbX}) the partial derivatives of the
instantaneous screw coordinates are%
\begin{equation}
\frac{\partial }{\partial q_{k}}\widehat{\mathbf{J}}{_{i,j}^{\text{b}}}=%
\frac{\partial }{\partial q_{k}}(\mathbf{C}_{i}^{-1}\mathbf{C}_{j})\mathbf{A}%
_{j}^{-1}\widehat{\mathbf{Y}}_{j}\mathbf{A}_{j}\mathbf{C}_{i}^{-1}\mathbf{C}%
_{j}+\mathbf{C}_{i}^{-1}\mathbf{C}_{j}\mathbf{A}_{j}^{-1}\widehat{\mathbf{Y}}%
_{j}\mathbf{A}_{j}\frac{\partial }{\partial q_{k}}(\mathbf{C}_{j}^{-1}%
\mathbf{C}_{i}).  \label{derTemp}
\end{equation}%
This can be evaluated with help of the POE formula (\ref{POEY}) as%
\begin{eqnarray}
\frac{\partial }{\partial q_{k}}(\mathbf{C}_{i}^{-1}\mathbf{C}_{j}) &=&\frac{%
\partial }{\partial q_{k}}(\mathbf{A}_{i}^{-1}\exp (-\widehat{\mathbf{Y}}%
_{i}q_{i})\cdots \exp (-\widehat{\mathbf{Y}}_{j+1}q_{j+1})\mathbf{A}_{j})) 
\notag \\
&=&-\mathbf{A}_{i}^{-1}\exp (-\widehat{\mathbf{Y}}_{i}q_{i})\cdots \exp (-%
\widehat{\mathbf{Y}}_{k+1}q_{k+1})\widehat{\mathbf{Y}}_{k}\exp (-\widehat{%
\mathbf{Y}}_{k}q_{k})\cdots \exp (-\widehat{\mathbf{Y}}_{j+1}q_{j+1})\mathbf{%
A}_{j}  \notag \\
&=&-\mathbf{C}_{i}^{-1}\mathbf{C}_{k}\mathbf{A}_{k}^{-1}\widehat{\mathbf{Y}}%
_{k}\mathbf{A}_{k}\mathbf{C}_{k}^{-1}\mathbf{C}_{j}=-\mathbf{C}_{i}^{-1}%
\mathbf{C}_{k}\mathbf{A}_{k}^{-1}\widehat{\mathbf{Y}}_{k}\mathbf{A}_{k}%
\mathbf{C}_{k}^{-1}\mathbf{C}_{i}\mathbf{C}_{i}^{-1}\mathbf{C}_{j}  \notag \\
&=&-\widehat{\mathbf{J}}{_{i,k}^{\text{b}}}\mathbf{C}_{i}^{-1}\mathbf{C}%
_{j},\ j\leq k\leq i,
\end{eqnarray}%
and in the same way follows that%
\begin{equation}
\frac{\partial }{\partial q_{k}}(\mathbf{C}_{j}^{-1}\mathbf{C}_{i})=\mathbf{C%
}_{j}^{-1}\mathbf{C}_{i}\widehat{\mathbf{J}}{_{i,j}^{\text{b}}},\ j\leq
k\leq i.
\end{equation}%
Inserted into (\ref{derTemp}) yields $\frac{\partial }{\partial q_{k}}%
\widehat{\mathbf{J}}{_{i,j}^{\text{b}}}=\widehat{\mathbf{J}}{_{i,j}^{\text{b}%
}}\widehat{\mathbf{J}}{_{i,k}^{\text{b}}}-\widehat{\mathbf{J}}{_{i,k}^{\text{%
b}}}\widehat{\mathbf{J}}{_{i,j}^{\text{b}}}$, and noting (\ref{LieBracketSE3}%
), the final expression is%
\begin{equation}
\begin{tabular}{|lll|}
\hline
&  &  \\ 
& $\dfrac{\partial \mathbf{J}{_{i,j}^{\text{b}}}}{\partial q_{k}}=[\mathbf{J}%
{_{i,j}^{\text{b}},\mathbf{J}_{i,k}^{\text{b}}}],~j<k\leq i.$ &  \\ 
&  &  \\ \hline
\end{tabular}
\label{der1}
\end{equation}%
Hence the partial derivative of the instantaneous joint screw $\mathbf{J}{%
_{i,j}^{\text{b}}}$ w.r.t. to $q_{k}$ is simply the screw product (\ref%
{ScrewProd}) of $\mathbf{J}{_{i,j}^{\text{b}}}$ and $\mathbf{J}{_{i,k}^{%
\text{b}}}$. The final expression for the acceleration attains a very
compact form%
\begin{equation}
\dot{\mathbf{V}}_{i}^{\text{b}}=\sum_{j\leq i}\mathbf{J}{_{i,j}^{\text{b}}}%
\ddot{q}_{j}+\sum_{j<k\leq i}[\mathbf{J}{_{i,j}^{\text{b}},\mathbf{J}{%
_{i,k}^{\text{b}}}]}\dot{q}_{j}\dot{q}_{k}.  \label{Vbdot}
\end{equation}%
Indeed the same result would be obtained using (\ref{JbX}) in terms of $%
\mathbf{Y}_{i}$. This expression has been derived, using different
notations, for instance in \cite{Brockett1984,MUBOLieGroup,Murray,Park1994}.

The equations (\ref{Vbdot}) can be summarized for all bodies $i=1,\ldots ,n$
using the system twist (\ref{VbAX}) and system Jacobian (\ref{JbSys}). To
this end, the derivative (\ref{der1}) is rewritten as 
\begin{equation}
\frac{\partial \mathbf{J}{_{i,j}^{\text{b}}}}{\partial q_{k}}=[\mathbf{J}{%
_{i,j}^{\text{b}},\mathbf{J}_{i,k}^{\text{b}}}]=\mathbf{Ad}_{\mathbf{C}%
_{i,k}}[\mathbf{J}{_{k,j}^{\text{b}},{^{k}}\mathbf{X}_{k}}]=-\mathbf{Ad}_{%
\mathbf{C}_{i,k}}\mathbf{ad}_{{{^{k}}\mathbf{X}_{k}}}\mathbf{J}{_{k,j}^{%
\text{b}}},~j<k\leq i  \label{dJbij}
\end{equation}%
so that 
\begin{equation*}
\dot{\mathbf{J}}{_{i,j}^{\text{b}}}=\sum_{j<k\leq i}[\mathbf{J}{_{i,j}^{%
\text{b}},\mathbf{J}_{i,k}^{\text{b}}}]\dot{q}_{k}=-\sum_{j<k\leq i}\mathbf{%
Ad}_{\mathbf{C}_{i,k}}\mathbf{ad}_{{{^{k}}\mathbf{X}_{k}}}\mathbf{J}{_{k,j}^{%
\text{b}}}\dot{q}_{k}.
\end{equation*}%
Noticing that $\mathbf{ad}_{{{^{k}}\mathbf{X}_{k}}}\mathbf{J}{_{k,k}^{\text{b%
}}}=\mathbf{0}$ the time derivative of the body-fixed system Jacobian
factors as%
\begin{equation}
\dot{\mathsf{J}}{^{\text{b}}}\left( \mathbf{q},\dot{\mathbf{q}}\right) =-%
\mathsf{A}^{\text{b}}\left( \mathbf{q}\right) \mathsf{a}^{\text{b}}\left( 
\dot{\mathbf{q}}\right) \mathsf{A}^{\text{b}}\left( \mathbf{q}\right) 
\mathsf{X}^{\text{b}}=-\mathsf{A}^{\text{b}}\left( \mathbf{q}\right) \mathsf{%
a}^{\text{b}}\left( \dot{\mathbf{q}}\right) \mathsf{J}^{\text{b}}\left( 
\mathbf{q}\right)
\end{equation}%
with $\mathsf{A}^{\text{b}}$ defined in (24) of \cite{Part1} and with%
\begin{equation}
\mathsf{a}^{\text{b}}\left( \dot{\mathbf{q}}\right) :=\mathrm{diag}~(\dot{q}%
_{1}\mathbf{ad}_{{{^{1}}\mathbf{X}_{1}}},\ldots ,\dot{q}_{n}\mathbf{ad}_{{{%
^{n}}\mathbf{X}_{n}}}).  \label{ab}
\end{equation}%
Hence the system acceleration is given in compact matrix form as%
\begin{equation}
\dot{\mathsf{V}}^{\text{b}}=\mathsf{J}{^{\text{b}}}\ddot{\mathbf{q}}-\mathsf{%
A}^{\text{b}}\mathsf{a}^{\text{b}}\mathsf{J}^{\text{b}}\dot{\mathbf{q}}=%
\mathsf{J}{^{\text{b}}}\ddot{\mathbf{q}}-\mathsf{A}^{\text{b}}\mathsf{a}^{%
\text{b}}\mathsf{V}^{\text{b}}.  \label{VbdotMat}
\end{equation}

\begin{remark}[Overall inverse kinematics solution]
The relation (\ref{VbdotMat}) gives rise to a solution of the inverse
kinematics problem on acceleration level, i.e. the generalized accelerations
for given configurations, twists, and accelerations of the bodies. The
unique solution is 
\begin{equation}
\ddot{\mathbf{q}}=((\mathsf{X}^{\text{b}})^{T}\mathsf{X}^{\text{b}})^{-1}(%
\mathsf{X}^{\text{b}})^{T}((\mathbf{I}-\mathsf{D}^{\text{b}})\dot{\mathsf{V}}%
^{\text{b}}+\mathsf{a}^{\text{b}}\mathsf{V}^{\text{b}})  \label{solqddot1}
\end{equation}%
which is indeed the time derivative of (26) in \cite{Part1}. In components
this gives the acceleration of the individual joints as $\ddot{q}_{i}={^{i}%
\mathbf{X}}_{i}^{T}(\dot{\mathbf{V}}_{i}^{\text{b}}-\mathbf{Ad}_{\mathbf{C}%
_{i,i-1}}\dot{\mathbf{V}}_{i-1}^{\text{b}}+\dot{q}_{i}[{^{i}\mathbf{X}}_{i},%
\mathbf{V}_{i}^{\text{b}}])/\left\Vert {^{i}\mathbf{X}}_{i}\right\Vert ^{2}$.
\end{remark}

A further time derivative of the twist yields the jerk of a body, which
requires a further partial derivative of the Jacobian. Starting from (\ref%
{der1}), and using the Jacobi identity (\ref{JacIdent}) and the bilinearity $%
\frac{\partial }{\partial q_{k}}[\mathbf{J}{_{i,j}^{\text{b}},\mathbf{J}{%
_{i,k}^{\text{b}}}]=}[\frac{\partial }{\partial q_{k}}\mathbf{J}{_{i,j}^{%
\text{b}},\mathbf{J}{_{i,k}^{\text{b}}}]+}[\mathbf{J}{_{i,j}^{\text{b}},%
\frac{\partial }{\partial q_{k}}\mathbf{J}{_{i,k}^{\text{b}}}]}$, the
non-zero second partial derivative is found as 
\begin{equation}
\frac{\partial ^{2}\mathbf{J}_{i,j}^{\text{b}}}{\partial q_{k}\partial q_{r}}%
=\left\{ 
\begin{array}{cl}
\lbrack \lbrack \mathbf{J}{_{i,j}^{\text{b}},\mathbf{J}{_{i,k}^{\text{b}}}],}%
\mathbf{J}{_{i,r}^{\text{b}}}], & j<k\leq r\leq i \\ 
\lbrack \lbrack \mathbf{J}{_{i,j}^{\text{b}},\mathbf{J}{_{i,r}^{\text{b}}}],}%
\mathbf{J}{_{i,k}^{\text{b}}}], & j<r<k\leq i%
\end{array}%
\right. .  \label{der2Jb}
\end{equation}%
This gives rise to an explicit form for the body-fixed jerk%
\begin{eqnarray}
\ddot{\mathbf{V}}_{i}^{\text{b}} &=&\sum_{j\leq i}\mathbf{J}{_{i,j}^{\text{b}%
}}\dddot{q}_{j}+2%
\hspace{-1ex}%
\sum_{j<k\leq i}[\mathbf{J}{_{i,j}^{\text{b}},\mathbf{J}{_{i,k}^{\text{b}}}]}%
\ddot{q}_{j}\dot{q}_{k}  \label{Vbddot} \\
&&+\sum_{j<k\leq i}[\mathbf{J}{_{i,j}^{\text{b}},\mathbf{J}{_{i,k}^{\text{b}}%
}]}\dot{q}_{j}\ddot{q}_{k}+2%
\hspace{-2ex}%
\sum_{j<k\leq r\leq i}[[\mathbf{J}{_{i,j}^{\text{b}},\mathbf{J}{_{i,k}^{%
\text{b}}}],}\mathbf{J}{_{i,r}^{\text{b}}}]\dot{q}_{j}\dot{q}_{k}\dot{q}_{r}.
\notag
\end{eqnarray}%
Thus only computationally simple nested screw products are required to
compute the terms that are quadratic and cubic in $\ddot{q}_{j}$, $\dot{q}%
_{k}$. The same applies to higher derivatives (that for instance are
required for motion planning). The explicit form of the $\nu $-th order
partial derivative was presented in \cite{ARK2014}%
\begin{equation}
\frac{\partial ^{\nu }\mathbf{J}_{i,j}^{\text{b}}}{\partial q_{\alpha
_{1}}\partial q_{\alpha _{2}}\cdots \partial q_{\alpha _{\nu }}}=[\ldots
\lbrack \lbrack \lbrack \mathbf{J}_{i,j}^{\text{b}}{,}\mathbf{J}_{i,\beta
_{1}}^{\text{b}}],\mathbf{J}_{i,\beta _{2}}^{\text{b}}{]},\mathbf{J}%
_{i,\beta _{3}}^{\text{b}}{]}\ldots ,\mathbf{J}_{i,\beta _{\nu }}^{\text{b}%
}],\ j<\beta _{1}\leq \beta _{2}\leq \cdots \leq \beta _{\nu }\leq i
\label{nthorderBody}
\end{equation}%
where $\beta _{1}\leq \beta _{2}\leq \cdots \leq \beta _{\nu }$ is the
ordered sequence of the indices $\alpha _{1},\ldots ,\alpha _{\nu }$.
Clearly the closed form expressions become very involved. Their explicit
determination can be avoided by recursive evaluation \cite{ARK2014}.

\subsection{Spatial Representation%
\label{secSpatialAcc}%
}

Proceeding in the same way as for (\ref{derTemp}) the partial derivative of
the spatial Jacobian is obtained as%
\begin{equation}
\begin{tabular}{|lll|}
\hline
&  &  \\ 
& $\dfrac{\partial \mathbf{J}{_{j}^{\text{s}}}}{\partial q_{k}}=[\mathbf{J}{%
_{k}^{\text{s}},\mathbf{J}_{j}^{\text{s}}}],~k<j.$ &  \\ 
&  &  \\ \hline
\end{tabular}
\label{Jsder}
\end{equation}%
Since the spatial representation $\mathbf{J}{_{j}^{\text{s}}}$ is intrinsic
to the joint $j$, rather than related to a body as is (\ref{JbX}), the time
derivative can be expressed as%
\begin{equation}
\dot{\mathbf{J}}{_{j}^{\text{s}}}=\sum_{k\leq j}[\mathbf{J}{_{k}^{\text{s}},%
\mathbf{J}_{j}^{\text{s}}}]\dot{q}^{k}=[\sum_{k\leq j}\mathbf{J}{_{k}^{\text{%
s}},\mathbf{J}_{j}^{\text{s}}}]\dot{q}_{k}=[\mathbf{V}{_{j}^{\text{s}},%
\mathbf{J}_{j}^{\text{s}}}].  \label{timeDerJs}
\end{equation}%
This relation reconfirms the special properties of spatial twists that are
advantageous for recursive implementations. It may be considered as an
extension of Euler's formula for the time derivative of vectors resolved in
moving frames to screws. For this reason Featherstone \cite%
{Featherstone1987,Featherstone2008} termed the Lie bracket the 'spatial
cross product'. The spatial acceleration is therewith%
\begin{equation}
\dot{\mathbf{V}}_{i}^{\text{s}}=\sum_{j\leq i}\mathbf{J}{_{j}^{\text{s}}}%
\ddot{q}_{j}+\sum_{k<j\leq i}[\mathbf{J}{_{k}^{\text{s}},\mathbf{J}{_{j}^{%
\text{s}}}]}\dot{q}_{j}\dot{q}_{k}=\sum_{j\leq i}\left( \mathbf{J}{_{j}^{%
\text{s}}}\ddot{q}_{j}+[\mathbf{V}{_{j}^{\text{s}},\mathbf{J}_{j}^{\text{s}}}%
]\dot{q}_{j}\right) .  \label{Vsdottemp}
\end{equation}%
In matrix form the overall spatial acceleration can be summarized as%
\begin{equation}
\dot{\mathsf{V}}^{\text{s}}=\mathsf{J}{^{\text{s}}}\ddot{\mathbf{q}}+{%
\mathsf{Lb}^{\text{s}}\mathrm{diag}~(\mathbf{J}_{1}^{\text{s}%
},\ldots ,\mathbf{J}_{n}^{\text{s}})\dot{\mathbf{q}}}  \label{VsdotMat}
\end{equation}%
with 
\begin{equation}
\mathsf{b}^{\text{s}}\left( \mathsf{V}^{\text{s}}\right) :=\mathrm{diag}~(%
\mathbf{ad}_{\mathbf{V}_{1}^{\text{s}}},\ldots ,\mathbf{ad}_{\mathbf{V}_{n}^{%
\text{s}}})
\end{equation}%
and $\mathsf{L}$ being the lower triangular block identity matrix. A
solution for $\ddot{\mathbf{q}}$ similar to (\ref{VbdotMat}) exists.

The second partial derivative of the spatial Jacobian is%
\begin{equation}
\frac{\partial ^{2}\mathbf{J}{_{i}^{\text{s}}}}{\partial q_{k}q_{j}}=\left\{ 
\begin{array}{cl}
\lbrack \mathbf{J}{_{k}^{\text{s}},}[\mathbf{J}{_{j}^{\text{s}},\mathbf{J}{%
_{i}^{\text{s}}}]}], & k<j<i \\ 
\lbrack \mathbf{J}{_{j}^{\text{s}},}[\mathbf{J}{_{k}^{\text{s}},\mathbf{J}{%
_{i}^{\text{s}}}]}], & j\leq k<i%
\end{array}%
\right. .  \label{Jsder2}
\end{equation}%
Therewith the spatial representation of the jerk of body $i$ is found as%
\begin{eqnarray}
\ddot{\mathbf{V}}_{i}^{\text{s}} &=&\sum_{j\leq i}%
\Big%
(\mathbf{J}{_{j}^{\text{s}}}\dddot{q}_{j}+2[\mathbf{V}{_{j}^{\text{s}}},%
\mathbf{J}{_{j}^{\text{s}}]}\ddot{q}_{j}+\sum_{k\leq j}[{\mathbf{J}{_{k}^{%
\text{s}}}}\ddot{q}_{k},\mathbf{J}{_{j}^{\text{s}}]}\dot{q}^{j}+[\mathbf{V}{%
_{j-1}^{\text{s}}}+\mathbf{V}{_{j}^{\text{s}}-\mathbf{V}{_{i}^{\text{s}}},[}%
\mathbf{V}{_{j}^{\text{s}}},\mathbf{J}{_{j}^{\text{s}}}]]\dot{q}_{j}%
\Big%
) \\
&=&\sum_{j\leq i}%
\Big%
(\mathbf{J}{_{j}^{\text{s}}}\dddot{q}_{j}+[[\mathbf{V}{_{j}^{\text{s}}},{%
\mathbf{J}{_{j}^{\text{s}}]}},\mathbf{V}{_{i}^{\text{s}}}-2\mathbf{V}{_{j}^{%
\text{s}}}]\dot{q}_{j}+[\sum_{k\leq j}{\mathbf{J}{_{k}^{\text{s}}}}\ddot{q}%
_{k}+[\mathbf{V}{_{j}^{\text{s}}},{\mathbf{J}{_{j}^{\text{s}}]}}\dot{q}_{j},{%
\mathbf{J}{_{j}^{\text{s}}}}]\dot{q}_{j}+2[\mathbf{V}{_{j}^{\text{s}}},{%
\mathbf{J}{_{j}^{\text{s}}}}]\ddot{q}_{j}%
\Big%
)
\end{eqnarray}%
The instantaneous joint screws (\ref{JsX}), and thus their derivatives (\ref%
{Jsder}) and (\ref{Jsder2}), are independent of a particular body. The
closed form of the $\nu $-th order partial derivative has been reported in 
\cite{MMTHighDer}%
\begin{eqnarray}
\frac{\partial ^{\nu }\mathbf{J}{_{i}^{\text{s}}}}{\partial q_{\alpha
_{1}}\partial q_{\alpha _{2}}\cdots \partial q_{\alpha _{\nu }}} &=&[\mathbf{%
J}_{\beta _{\nu }}^{\text{s}},[\mathbf{J}_{\beta _{\nu -1}}^{\text{s}},[%
\mathbf{J}_{\beta _{\nu -2}}^{\text{s}},\ldots \lbrack \mathbf{J}_{\beta
_{1}}^{\text{s}},\mathbf{J}_{i}^{\text{s}}]\ldots ]]],\ \beta _{\nu }\leq
\beta _{\nu -1}\leq \cdots \leq \beta _{1}<i%
\vspace{1ex}
\label{nthorder} \\
&=&\mathbf{ad}_{\mathbf{J}_{\beta _{\nu }}^{\text{s}}}\mathbf{ad}_{\mathbf{J}%
_{\beta _{\nu -1}}^{\text{s}}}\mathbf{ad}_{\mathbf{J}_{\beta _{\nu -2}}^{%
\text{s}}}\cdots \mathbf{ad}_{\mathbf{J}_{\beta _{1}}^{\text{s}}}\mathbf{J}%
_{i}^{\text{s}},\ \beta _{\nu }\leq \beta _{\nu -1}\leq \cdots \leq \beta
_{1}<i%
\vspace{1ex}
\notag \\
&=&[\mathbf{J}_{\beta _{\nu }}^{\text{s}},\frac{\partial ^{\nu -1}\mathbf{J}%
_{i}^{\text{s}}}{\partial q_{\beta _{1}}\partial q_{\beta _{2}}\cdots
\partial q_{\beta _{\nu -1}}}],\ \beta _{\nu }\leq \beta _{\nu -1}<i  \notag
\end{eqnarray}%
where again $\beta _{\nu }\leq \beta _{\nu -1}\leq \cdots \leq \beta _{1}$
is the ordered sequence of the indices $\alpha _{1},\ldots ,\alpha _{\nu }$.
The last form in (\ref{nthorder}) allows for a recursive determination.
Moreover, a recursive formulation for the time derivative of spatial twists
has been reported in \cite{ConstraintsMMT2016}. Together with the very
concise form (\ref{timeDerJs}) this makes the spatial representation
computationally very attractive.

\subsection{Hybrid Form}

The results in section \ref{secDerBody} can be carried over to the hybrid
twist making use of the relation (\ref{JhJb}). As in (\ref{Vw}), denote with 
${{\overset{v}{\mathbf{J}}}{_{i,k}^{\text{h}}}}$ and ${{\overset{\omega }{%
\mathbf{J}}}{_{i,k}^{\text{h}}}}$ the screw coordinate vectors comprising
respectively the linear and angular part of the column of the hybrid
Jacobian so that ${\mathbf{J}_{i,k}^{\text{h}}}={{\overset{\omega }{\mathbf{J%
}}}{_{i,k}^{\text{h}}+\overset{v}{\mathbf{J}}}{_{i,k}^{\text{h}}}}$. Then%
\begin{eqnarray}
\frac{\partial {\mathbf{J}_{i,j}^{\text{h}}}}{\partial q_{k}} &=&\frac{%
\partial \mathbf{Ad}_{\mathbf{R}_{i}}}{\partial q_{k}}\mathbf{J}_{i,j}^{%
\text{b}}+\mathbf{Ad}_{\mathbf{R}_{i}}\frac{\partial \mathbf{J}_{i,j}^{\text{%
b}}}{\partial q_{k}}=\mathbf{ad}_{{{\overset{\omega }{\mathbf{J}}}{_{i,k}^{%
\text{h}}}}}\mathbf{Ad}_{\mathbf{R}_{i}}\mathbf{J}_{i,j}^{\text{b}}+\mathbf{%
Ad}_{\mathbf{R}_{i}}[\mathbf{J}_{i,j}^{\text{b}},\mathbf{J}_{i,k}^{\text{b}}]
\notag \\
&=&[{{\overset{\omega }{\mathbf{J}}}{_{i,k}^{\text{h}}}}{,\mathbf{J}_{i,j}^{%
\text{h}}]+[\mathbf{J}_{i,j}^{\text{h}},\mathbf{J}_{i,k}^{\text{h}}]=[%
\mathbf{J}_{i,j}^{\text{h}},\mathbf{J}_{i,k}^{\text{h}}-{\overset{\omega }{%
\mathbf{J}}}{_{i,k}^{\text{h}}}]}
\end{eqnarray}%
and thus%
\begin{equation}
\begin{tabular}{|lll|}
\hline
&  &  \\ 
& $\dfrac{\partial {\mathbf{J}_{i,j}^{\text{h}}}}{\partial q_{k}}{=}{[%
\mathbf{J}_{i,j}^{\text{h}},{\overset{v}{\mathbf{J}}}{_{i,k}^{\text{h}}}]}=-%
\mathbf{ad}_{{{\overset{v}{\mathbf{J}}}{_{i,k}^{\text{h}}}}}{\mathbf{J}%
_{ij}^{\text{h}},\ \ j\leq k\leq i.}$ &  \\ 
&  &  \\ \hline
\end{tabular}
\label{Jhder}
\end{equation}%
The similarity to (\ref{der1}) is apparent. The difference is that the
convective term due to the angular motion is missing, which is why only $%
\overset{v}{\mathbf{J}}$ appears. The time derivative of the hybrid Jacobian
can thus be expressed as%
\begin{equation}
\dot{\mathbf{J}}{_{i,j}^{\text{h}}}=\sum_{k\leq j}{[\mathbf{J}_{i,j}^{\text{h%
}},{\overset{v}{\mathbf{J}}}{_{i,k}^{\text{h}}}]}\dot{q}_{k}={[\mathbf{J}%
_{i,j}^{\text{h}},\Delta {\overset{v}{\mathbf{V}}}{_{j-1,i}^{\text{h}}}]}%
\text{ }
\end{equation}%
where $\Delta {\mathbf{V}{_{j-1,i}^{\text{h}}}}:={\mathbf{V}{_{i}^{\text{h}}}%
}-\mathbf{Ad}_{\mathbf{r}_{i,j-1}}{\mathbf{V}{_{j-1}^{\text{h}}}}$ is the
relative hybrid twist of body $j-1$ and $i$ as observed in the BFR on body $i
$. A simpler relation is obtained by directly differentiating (\ref{Jh})%
\begin{eqnarray}
\dot{\mathbf{J}}{_{i,j}^{\text{h}}} &=&(\mathbf{ad}_{\dot{\mathbf{r}}_{i,j}}+%
\mathbf{Ad}_{\mathbf{r}_{i,j-1}}\mathbf{ad}_{\mathbf{\omega }_{j}^{\text{s}%
}}){^{0}}\mathbf{X}_{j}^{j}  \label{Jhdot} \\
&=&\mathbf{Ad}_{\mathbf{r}_{i,j-1}}(\mathbf{ad}_{\dot{\mathbf{d}}_{i,j}}+%
\mathbf{ad}_{\mathbf{\omega }_{j}^{\text{s}}}){^{0}}\mathbf{X}_{j}^{j}=%
\mathbf{Ad}_{\mathbf{r}_{i,j-1}}(\mathbf{ad}_{\mathbf{V}_{j}^{\text{s}}}-%
\mathbf{ad}_{\dot{\mathbf{r}}_{i}}){^{0}}\mathbf{X}_{j}^{j}.  \notag
\end{eqnarray}%
This yields the following explicit expressions for the hybrid acceleration%
\begin{align}
\dot{\mathbf{V}}_{i}^{\text{h}}=& \sum_{j\leq i}\mathbf{J}{_{{i,}j}^{\text{h}%
}}\ddot{q}_{j}+\sum_{j\leq k\leq i}{[\mathbf{J}_{i,j}^{\text{h}},{\overset{v}%
{\mathbf{J}}}{_{i,k}^{\text{h}}}]}\dot{q}_{j}\dot{q}_{k}=\sum_{j\leq i}(%
\mathbf{J}{_{{i,}j}^{\text{h}}}\ddot{q}_{j}+{[\mathbf{J}_{i,j}^{\text{h}%
},\Delta {\overset{v}{\mathbf{V}}}{_{j-1,i}^{\text{h}}}]}\dot{q}_{j}) \\
=& \sum_{j\leq i}(\mathbf{J}{_{{i,}j}^{\text{h}}}\ddot{q}_{j}+(\mathbf{ad}_{%
\dot{\mathbf{r}}_{i,j}}+\mathbf{Ad}_{\mathbf{r}_{i,j}}\mathbf{ad}_{\mathbf{%
\omega }_{j}^{\text{s}}}){^{0}}\mathbf{X}_{j}^{j}\dot{q}_{j}).
\end{align}%
For the second derivative it is simplest to start from (\ref{Jhdot}), and a
straightforward calculation yields%
\begin{equation}
\ddot{\mathbf{J}}{_{i,j}^{\text{h}}}=\left( \mathbf{ad}_{\ddot{\mathbf{r}}%
_{i,j}}+2\mathbf{ad}_{\dot{\mathbf{r}}_{i,j}}\mathbf{ad}_{\mathbf{\omega }%
_{j}^{\text{s}}}+\mathbf{Ad}_{\mathbf{r}_{i,j}}(\mathbf{ad}_{\dot{\mathbf{%
\omega }}_{j}^{\text{s}}}+\mathbf{ad}_{\mathbf{\omega }_{j}^{\text{s}}}%
\mathbf{ad}_{\mathbf{\omega }_{j}^{\text{s}}})\right) {^{0}}\mathbf{X}%
_{j}^{j}.
\end{equation}%
The jerk in hybrid representation can thus be written as%
\begin{align}
\ddot{\mathbf{V}}_{i}^{\text{h}}=& \sum_{j\leq i}\left( \mathbf{J}{_{{i,}j}^{%
\text{h}}}\dddot{q}_{j}+2\mathbf{ad}_{\dot{\mathbf{r}}_{i,j}}\ddot{q}_{j}+%
\big%
(\mathbf{ad}_{\ddot{\mathbf{r}}_{i,j}}+2\mathbf{ad}_{\dot{\mathbf{r}}_{i,j}}%
\mathbf{ad}_{\mathbf{\omega }_{j}^{\text{s}}}%
\big%
)\dot{q}_{j}\right.  \\
& \ \ \ \ \ \ +\left. \mathbf{Ad}_{\mathbf{r}_{i,j}}%
\big%
(2\mathbf{ad}_{\mathbf{\omega }_{j}^{\text{s}}}\ddot{q}_{j}+\mathbf{ad}_{%
\dot{\mathbf{\omega }}_{j}^{\text{s}}}+\mathbf{ad}_{\mathbf{\omega }_{j}^{%
\text{s}}}\mathbf{ad}_{\mathbf{\omega }_{j}^{\text{s}}}%
\big%
)\dot{q}_{j}){^{0}}\mathbf{X}_{j}^{j}\right) .
\end{align}%
These are the core relation in the so-called 'spatial vector' formulation
(i.e. using the hybrid representation of twists) \cite%
{Jian1991,Jian1995,Fijany1995,LillyOrin1991,Rodriguez1987,Rodriguez1992}. In
this context the Lie bracket, respectively screw product, (\ref{ScrewProd})
has been termed the 'spatial cross product' \cite%
{Featherstone1987,Featherstone2008}.

\subsection{Mixed Representation}

With (\ref{Vmbhs}), employing the results for the mixed representation,
yields

\begin{equation}
\dot{\mathbf{J}}{_{ij}^{\text{m}}}=\left( 
\begin{array}{cc}
\mathbf{R}_{i}^{T} & \mathbf{0} \\ 
\mathbf{0} & \mathbf{I}%
\end{array}%
\right) \dot{\mathbf{J}}{_{ij}^{\text{h}},\ \ }\dot{\mathbf{V}}_{i}^{\text{m}%
}=\left( 
\begin{array}{cc}
\mathbf{R}_{i}^{T} & \mathbf{0} \\ 
\mathbf{0} & \mathbf{I}%
\end{array}%
\right) \dot{\mathbf{V}}_{i}^{\text{h}},\ \ \ddot{\mathbf{V}}_{i}^{\text{m}%
}=\left( 
\begin{array}{cc}
\mathbf{R}_{i}^{T} & \mathbf{0} \\ 
\mathbf{0} & \mathbf{I}%
\end{array}%
\right) \ddot{\mathbf{V}}_{i}^{\text{h}}.
\end{equation}

\section{Newton-Euler Equations in Various Representations%
\label{secNE}%
}

\subsection{Spatial Representation%
\label{secNESpatial}%
}

Consider a rigid body with body-fixed BFR $\mathcal{F}_{\text{b}}=\{\Omega ;%
\vec{e}_{\text{b},1},\vec{e}_{\text{b},2},\vec{e}_{\text{b},3}\}$ located at
an arbitrary point $\Omega $. Denote the inertia matrix w.r.t. this BFR with 
$\mathbf{M}^{\text{b}}$, ref. (\ref{Mb}). The configuration of the BFR $%
\mathcal{F}_{\text{b}}$ is described by $\mathbf{C}=\left( \mathbf{R},%
\mathbf{r}\right) $. The spatial inertia matrix expressed in the IFR is then%
\begin{equation}
\mathbf{M}^{\text{s}}=\mathbf{Ad}_{\mathbf{C}}^{-T}\mathbf{M}^{\text{b}}%
\mathbf{Ad}_{\mathbf{C}}^{-1}.  \label{Ms}
\end{equation}%
The spatial canonical momentum co-screw $\mathbf{\Pi }^{\text{s}}=\left( 
\mathbf{L}^{\text{s}},\mathbf{P}^{\text{s}}\right) ^{T}\in se^{\ast }\left(
3\right) $, conjugate to the spatial twist, is thus%
\begin{equation}
\mathbf{\Pi }^{\text{s}}=\mathbf{M}^{\text{s}}\mathbf{V}^{\text{s}}=\mathbf{%
Ad}_{\mathbf{C}}^{-T}\mathbf{M}^{\text{b}}\mathbf{Ad}_{\mathbf{C}}^{-1}%
\mathbf{V}^{\text{s}}=\mathbf{Ad}_{\mathbf{C}}^{-T}\mathbf{\Pi }^{\text{b}}.
\label{Pis}
\end{equation}%
The momentum balance yields the Newton-Euler (NE) equations in spatial
representation, which attains the simple form%
\begin{equation}
\dot{\mathbf{\Pi }}^{\text{s}}=\mathbf{W}^{\text{s}}  \label{Pisdot}
\end{equation}%
where $\mathbf{W}^{\text{s}}=\left( \mathbf{t}^{\text{s}},\mathbf{f}^{\text{s%
}}\right) ^{T}$ is the applied wrench, with spatial torque $\mathbf{t}^{%
\text{s}}\equiv {^{0}}\mathbf{t}^{0}$ and force $\mathbf{f}^{\text{s}}\equiv 
{^{0}}\mathbf{f}$, both measured and resolved in the IFR. The momentum
balance equation (\ref{Pisdot}) is the simplest form possible, which is
achieved by using the spatial representation of twist, wrench, and momentum.
Firstly, it does not involve any vectorial operation, e.g. cross products.
Secondly, it is also numerically advantageous: any numerical discretization
of the ODE (\ref{Pisdot}) easily preserves the spatial momentum in the
absence of external wrenches. This has been discussed already by Borri in 
\cite{Bottasso1998}. In this context the spatial formulation is called the
fixed pole equation. In a recent paper \cite{GacesaJelenic2015} the
advantages of this form are exploited for geometrically exact modeling of
beams.

The explicit and compact form in terms of the spatial twist is found,
introducing (\ref{Pis}) and using 
\begin{equation}
\dot{\mathbf{M}}^{\text{s}}=-\mathbf{ad}_{\mathbf{V}^{\text{s}}}^{T}\mathbf{M%
}^{\text{s}}-\mathbf{M}^{\text{s}}\mathbf{ad}_{\mathbf{V}^{\text{s}}}
\label{Msdot}
\end{equation}%
along with $\mathbf{ad}_{\mathbf{V}^{\text{s}}}\mathbf{V}^{\text{s}}=\mathbf{%
0}$, as 
\begin{equation}
\begin{tabular}{|lll|}
\hline
&  &  \\ 
& $\mathbf{W}^{\text{s}}=\mathbf{M}^{\text{s}}\dot{\mathbf{V}}^{\text{s}}-%
\mathbf{ad}_{\mathbf{V}^{\text{s}}}^{T}\mathbf{M}^{\text{s}}\mathbf{V}^{%
\text{s}}.$ &  \\ 
&  &  \\ \hline
\end{tabular}
\label{NEs}
\end{equation}

\begin{remark}
Writing (\ref{NEs}) as $\mathbf{W}^{\text{s}}=\mathbf{M}^{\text{s}}\dot{%
\mathbf{V}}^{\text{s}}+\mathbf{C}^{\text{s}}\mathbf{V}^{\text{s}}$ (with $%
\mathbf{C}^{\text{s}}:=-\mathbf{ad}_{\mathbf{V}^{\text{s}}}^{T}\mathbf{M}^{%
\text{s}}$) shows that $\dot{\mathbf{M}}^{\text{s}}-2\mathbf{C}^{\text{s}}=%
\mathbf{ad}_{\mathbf{V}^{\text{s}}}^{T}\mathbf{M}^{\text{s}}-\mathbf{M}^{%
\text{s}}\mathbf{ad}_{\mathbf{V}^{\text{s}}}$ is skew symmetric. This
property is called the skew symmetry of the motion equations \cite{Murray}.
\end{remark}

\subsection{Body-fixed Representation%
\label{secNEBody}%
}

Let $\mathcal{F}_{\text{c}}=\{C;\vec{e}_{\text{c},1},\vec{e}_{\text{c},2},%
\vec{e}_{\text{c},3}\}$ be a body-fixed frame located at the COM. Its
configuration is described by $C_{\text{c}}=\left( \mathbf{R}_{\text{c}},%
\mathbf{r}_{\text{c}}\right) $. The body-fixed twist of the COM frame is
denoted $\mathbf{V}_{\text{c}}^{\text{b}}=\left( \mathbf{\omega }_{\text{c}%
}^{\text{b}},\mathbf{v}_{\text{c}}^{\text{b}}\right) ^{T}$ with $\widetilde{%
\mathbf{\omega }}_{\text{c}}^{\text{b}}=\mathbf{R}_{\text{c}}^{T}\dot{%
\mathbf{R}}_{\text{c}},\mathbf{v}_{\text{c}}^{\text{b}}=\mathbf{R}_{\text{c}%
}^{T}\dot{\mathbf{r}}_{\text{c}}$. The inertia matrix w.r.t. this COM frame
is denoted%
\begin{equation}
\mathbf{M}_{\text{c}}^{\text{b}}=\left( 
\begin{array}{cc}
\mathbf{\Theta }_{\text{c}} & \mathbf{0} \\ 
\mathbf{0} & m\mathbf{I}%
\end{array}%
\right)  \label{Mbc}
\end{equation}%
with the body mass $m$ and the inertia tensor $\mathbf{\Theta }_{\text{c}}$
expressed in the body-fixed COM frame $\mathcal{F}_{\text{c}}$. Let $S_{%
\text{bc}}=\left( \mathbf{R}_{\text{bc}},{^{\text{b}}}\mathbf{d}_{\text{bc}%
}\right) \in SE\left( 3\right) $ be the transformation from the COM frame $%
\mathcal{F}_{\text{c}}$ to the BFR $\mathcal{F}_{\text{b}}$. Here ${^{\text{b%
}}}\mathbf{d}_{\text{bc}}$ is the position vector from the BFR to the COM
resolved in the BFR. Then the configuration of $\mathcal{F}_{\text{c}}$ is
given in terms of that of the BFR by $\mathbf{C}_{\text{c}}=\mathbf{CS}_{%
\text{bc}}$. The inertia matrix w.r.t. to the general BFR $\mathcal{F}_{%
\text{b}}$ is%
\begin{eqnarray}
\mathbf{M}^{\text{b}} &=&\mathbf{Ad}_{\mathbf{S}_{\text{bc}}}^{-T}\mathbf{M}%
_{\text{c}}^{\text{b}}\mathbf{Ad}_{\mathbf{S}_{\text{bc}}}^{-1}  \notag \\
&=&\left( 
\begin{array}{cc}
\mathbf{\Theta }_{\text{b}} & m{^{\text{b}}\widetilde{\mathbf{d}}}_{\text{bc}%
} \\ 
-m{^{\text{b}}}\widetilde{\mathbf{d}}_{\text{bc}} & m\mathbf{I}%
\end{array}%
\right)  \label{Mb}
\end{eqnarray}%
with $\mathbf{\Theta }_{\text{b}}=\mathbf{R}_{\text{bc}}\mathbf{\Theta }_{%
\text{c}}\mathbf{R}_{\text{bc}}^{T}-m\widetilde{\mathbf{d}}_{\text{bc}}^{2}$
(which is the parallel axes theorem).

The momentum co-screw represented in the body-fixed RFR $\mathcal{F}_{\text{b%
}}$ is $\mathbf{\Pi }^{\text{b}}=\mathbf{M}^{\text{b}}\mathbf{V}^{\text{b}}$%
. The frame transformation of (\ref{NEs}) to the BFR $\mathcal{F}_{\text{b}}$
yields the body-fixed momentum balance represented in $\mathcal{F}_{\text{b}}
$ in the concise form%
\begin{equation}
\begin{tabular}{|llll|}
\hline
&  &  &  \\ 
& $\mathbf{W}^{\text{b}}%
\hspace{-2ex}%
$ & $=\dot{\mathbf{\Pi }}^{\text{b}}-\mathbf{ad}_{\mathbf{V}^{\text{b}}}^{T}%
\mathbf{\Pi }^{\text{b}}$ &  \\ 
&  & $=\mathbf{M}^{\text{b}}\dot{\mathbf{V}}^{\text{b}}-\mathbf{ad}_{\mathbf{%
V}^{\text{b}}}^{T}\mathbf{M}^{\text{b}}\mathbf{V}^{\text{b}}$ &  \\ 
&  &  &  \\ \hline
\end{tabular}
\label{NEb}
\end{equation}%
with the applied wrench $\mathbf{W}^{\text{b}}=\left( \mathbf{t}^{\text{b}},%
\mathbf{f}^{\text{b}}\right) ^{T}$ in body-fixed representation. The
equations (\ref{NEb}) are formally identical to the spatial equations (\ref%
{NEs}). Written separately, this yields the NE equations expressed in an
arbitrary body-fixed BFR%
\begin{eqnarray}
\mathbf{\Theta }_{\text{b}}\dot{\mathbf{\omega }}^{\text{b}}+\widetilde{%
\mathbf{\omega }}{^{\text{b}}}\mathbf{\Theta }_{\text{b}}\mathbf{\omega }^{%
\text{b}}-m{^{\text{b}}}\left( \dot{\mathbf{v}}^{\text{b}}+\widetilde{%
\mathbf{\omega }}{^{\text{b}}}\mathbf{v}{^{\text{b}}}\right) \widetilde{%
\mathbf{d}}_{\text{bc}} &=&\mathbf{t}^{\text{b}}  \label{NEbexp1} \\
m%
\big%
(\dot{\mathbf{v}}^{\text{b}}+\widetilde{\mathbf{\omega }}{^{\text{b}}}%
\mathbf{v}^{\text{b}}+(\dot{\widetilde{\mathbf{\omega }}}{^{\text{b}}}+%
\widetilde{\mathbf{\omega }}{^{\text{b}}}\widetilde{\mathbf{\omega }}{^{%
\text{b}})^{\text{b}}}\mathbf{d}_{\text{bc}}%
\big%
) &=&\mathbf{f}^{\text{b}}.  \label{NEbexp2}
\end{eqnarray}%
When using the COM frame as special case, the momentum represented in the
body-fixed COM frame is $\mathbf{\Pi }_{\text{c}}^{\text{b}}=\mathbf{M}_{%
\text{c}}^{\text{b}}\mathbf{V}_{\text{c}}^{\text{b}}$, and the momentum
balance yields%
\begin{equation}
\mathbf{W}_{\text{c}}^{\text{b}}=\mathbf{M}_{\text{c}}^{\text{b}}\dot{%
\mathbf{V}}_{\text{c}}^{\text{b}}-\mathbf{ad}_{\mathbf{\omega }_{\text{c}}^{%
\text{b}}}^{T}\mathbf{M}_{\text{c}}^{\text{b}}\mathbf{V}_{\text{c}}^{\text{b}%
}.  \label{NEb0}
\end{equation}%
Written in components, this yields the NE equations represented in the COM
frame%
\begin{eqnarray}
\mathbf{\Theta }_{\text{c}}\dot{\mathbf{\omega }}_{\text{c}}^{\text{b}}+%
\widetilde{\mathbf{\omega }}_{\text{c}}^{\text{b}}\mathbf{\Theta }_{\text{c}}%
\mathbf{\omega }_{\text{c}}^{\text{b}} &=&\mathbf{t}_{\text{c}}^{\text{b}}
\label{NEb0exp1} \\
m%
\big%
(\dot{\mathbf{v}}_{\text{c}}^{\text{b}}+\widetilde{\mathbf{\omega }}_{\text{c%
}}^{\text{b}}\mathbf{v}_{\text{c}}^{\text{b}}%
\big%
) &=&\mathbf{f}_{\text{c}}^{\text{b}}.  \label{NEb0exp2}
\end{eqnarray}%
Noticeably the angular and translational momentum equations are coupled even
though the COM is used as reference. This is due to using body-fixed twists.

\subsection{Hybrid Form}

The hybrid twist ${\mathbf{V}}{_{\text{c}}^{\text{h}}}=\left( \mathbf{\omega 
}^{\text{s}},\dot{\mathbf{r}}_{\text{c}}\right) ^{T}$ of the COM frame is
related to the body-fixed twist by $\mathbf{Ad}_{\mathbf{R}_{\text{c}}}^{-1}%
\mathbf{V}_{\text{c}}^{\text{b}}$, see (\ref{VhVb}), where $\mathbf{R}_{%
\text{c}}$ is the absolute rotation matrix of $\mathcal{F}_{\text{c}}$ in $%
C_{\text{c}}$. The hybrid momentum screw is thus $\mathbf{\Pi }_{\text{c}}^{%
\text{h}}=\mathbf{M}_{\text{c}}^{\text{h}}\mathbf{V}_{\text{c}}^{\text{h}}$,
where the hybrid representation of the inertia matrix is%
\begin{equation}
\mathbf{M}_{\text{c}}^{\text{h}}=\mathbf{Ad}_{\mathbf{R}_{\text{c}}}^{-T}%
\mathbf{M}_{\text{c}}^{\text{b}}\mathbf{Ad}_{\mathbf{R}_{\text{c}%
}}^{-1}=\left( 
\begin{array}{cc}
\mathbf{\Theta }_{\text{c}}^{\text{h}} & \mathbf{0} \\ 
\mathbf{0} & m\mathbf{I}%
\end{array}%
\right) ,\ \ \mathbf{\Theta }_{\text{c}}^{\text{h}}=\mathbf{R_{\text{c}}}%
\mathbf{\Theta }_{\text{c}}\mathbf{R}_{\text{c}}^{T}.  \label{M0h}
\end{equation}%
The hybrid momentum balance w.r.t. the COM follows from $\dot{\mathbf{\Pi }}%
_{\text{c}}^{\text{h}}=\mathbf{W}_{\text{c}}^{\text{h}}$. Using $\dot{%
\mathbf{M}}_{\text{c}}^{\text{h}}=-\mathbf{ad}_{\mathbf{\omega }^{\text{s}%
}}^{T}\mathbf{M}_{\text{c}}^{\text{h}}-\mathbf{M}_{\text{c}}^{\text{h}}%
\mathbf{ad}_{\mathbf{\omega }^{\text{s}}}$ yields%
\begin{equation}
\mathbf{W}_{\text{c}}^{\text{h}}=\mathbf{M}_{\text{c}}^{\text{h}}\dot{%
\mathbf{V}}_{\text{c}}^{\text{h}}+\mathbf{ad}_{\mathbf{\omega }^{\text{s}}}%
\mathbf{M}_{\text{c}}^{\text{h}}{\overset{\omega }{\mathbf{V}}}{_{\text{c}}^{%
\text{h}}}  \label{NEh0}
\end{equation}%
with ${\overset{\omega }{\mathbf{V}}}{_{\text{c}}^{\text{h}}}=\left( \mathbf{%
\omega }^{\text{s}},\mathbf{0}\right) ^{T}$ (notice $-\mathbf{ad}_{\mathbf{%
\omega }^{\text{s}}}^{T}=\mathbf{ad}_{\mathbf{\omega }^{\text{s}}}$).
Writing (\ref{NEh0}) separately for the angular and linear momentum balance 
\begin{eqnarray}
\mathbf{\Theta }_{\text{c}}^{\text{h}}\dot{\mathbf{\omega }}^{\text{s}}+%
\widetilde{\mathbf{\omega }}^{\text{s}}\mathbf{\Theta }_{\text{c}}^{\text{h}}%
\mathbf{\omega }^{\text{s}} &=&\mathbf{t}_{\text{c}}^{\text{h}}
\label{NEhCOM} \\
m\dot{\mathbf{r}}_{\text{c}} &=&\mathbf{f}_{\text{c}}^{\text{h}}
\label{NEhCOM2}
\end{eqnarray}%
shows that the hybrid NE equations w.r.t. the COM are indeed decoupled. Here 
$\mathbf{W}_{\text{c}}^{\text{h}}=(\mathbf{t}_{\text{c}}^{\text{h}},\mathbf{f%
}_{\text{c}}^{\text{h}})^{T}$ denotes the hybrid wrench measured in the COM
frame and resolved in the IFR.

Now consider the arbitrary body-fixed BFR $\mathcal{F}_{\text{b}}$ with
configuration $C=\left( \mathbf{R},\mathbf{r}\right) $. The hybrid twist $%
\mathbf{V}^{\text{h}}=\left( \mathbf{\omega }^{\text{s}},\dot{\mathbf{r}}%
\right) ^{T}$ measured at this RFR is $\mathbf{V}^{\text{h}}=\mathbf{Ad}_{%
\mathbf{d}_{\text{bc}}}\mathbf{V}_{\text{c}}^{\text{h}}$, with the
displacement vector $\mathbf{d}_{\text{bc}}$ from BFR to COM resolved in the
IFR. The hybrid mass matrix w.r.t. to the BFR $\mathcal{F}_{\text{b}}$ is
found as 
\begin{equation}
\mathbf{M}^{\text{h}}=\mathbf{Ad}_{\mathbf{d}_{\text{bc}}}^{-T}\mathbf{M}_{%
\text{c}}^{\text{h}}\mathbf{Ad}_{\mathbf{d}_{\text{bc}}}^{-1}=\left( 
\begin{array}{cc}
\mathbf{\Theta }^{\text{h}} & m\widetilde{\mathbf{d}}_{\text{bc}} \\ 
-m\widetilde{\mathbf{d}}_{\text{bc}} & m\mathbf{I}%
\end{array}%
\right) ,\ \ \mathbf{\Theta }^{\text{h}}=\mathbf{\Theta }_{\text{c}}^{\text{h%
}}-m\widetilde{\mathbf{d}}_{\text{bc}}^{2}.
\end{equation}%
The momentum balance in hybrid representation w.r.t. an arbitrary BFR%
\begin{equation}
\dot{\mathbf{\Pi }}^{\text{h}}=\mathbf{W}^{\text{h}}  \label{Pihdot}
\end{equation}%
is found, using $\dot{\mathbf{Ad}}_{\mathbf{d}_{\text{bc}}}^{-1}=-\mathbf{ad}%
_{\dot{\mathbf{d}}_{\text{bc}}}=\mathbf{Ad}_{\mathbf{d}_{\text{bc}}}^{-1}%
\mathbf{ad}_{\mathbf{\omega }^{\text{s}}}-\mathbf{ad}_{\mathbf{\omega }^{%
\text{s}}}\mathbf{Ad}_{\mathbf{d}_{\text{bc}}}^{-1}$ to evaluate (\ref%
{Pihdot}), as%
\begin{equation}
\begin{tabular}{|lll|}
\hline
&  &  \\ 
& $\mathbf{W}^{\text{h}}=\mathbf{M}^{\text{h}}\dot{\mathbf{V}}^{\text{h}}+%
\mathbf{ad}_{\mathbf{\omega }^{\text{s}}}\mathbf{M}^{\text{h}}{\overset{%
\omega }{\mathbf{V}}}{^{\text{h}}}.$ &  \\ 
&  &  \\ \hline
\end{tabular}
\label{NEh}
\end{equation}%
Separating the angular and translational part results in%
\begin{eqnarray}
\mathbf{\Theta }^{\text{h}}\dot{\mathbf{\omega }}^{\text{s}}+\widetilde{%
\mathbf{\omega }}^{\text{s}}\mathbf{\Theta }^{\text{h}}\mathbf{\omega }^{%
\text{s}}+m\widetilde{\mathbf{d}}_{\text{bc}}\ddot{\mathbf{r}} &=&\mathbf{t}%
^{\text{h}}  \label{NEhexp1} \\
m(\ddot{\mathbf{r}}+(\dot{\widetilde{\mathbf{\omega }}}{^{\text{s}}}+%
\widetilde{\mathbf{\omega }}^{\text{s}}\widetilde{\mathbf{\omega }}^{\text{s}%
})\mathbf{d}_{\text{bc}}) &=&\mathbf{f}^{\text{h}}.  \label{NEhexp2}
\end{eqnarray}%
These are simpler than the body-fixed equations (\ref{NEbexp1}) and (\ref%
{NEbexp2}). Finally notice that $\mathbf{f}^{\text{h}}=\mathbf{f}^{\text{s}}$%
.

\subsection{Mixed Form}

The mixed twists $\mathbf{V}^{\text{m}}=\left( \mathbf{\omega }^{\text{b}},%
\dot{\mathbf{r}}\right) ^{T}$ consists of the body-fixed angular velocity $%
\mathbf{\omega }^{\text{b}}$, i.e. measured and resolved in the BFR $%
\mathcal{F}_{\text{b}}$, and the translational velocity $\dot{\mathbf{r}}$
measured at the BFR $\mathcal{F}_{\text{b}}$ and resolved in the IFR. The NE
equations for the mixed representation w.r.t. a general BFR are directly
found by combining (\ref{NEbexp1}) and (\ref{NEhexp2}), with $\ddot{\mathbf{r%
}}=\dot{\mathbf{v}}^{\text{b}}+\widetilde{\mathbf{\omega }}{^{\text{b}}}%
\mathbf{v}{^{\text{b}}}$, 
\begin{eqnarray}
\mathbf{\Theta }^{\text{b}}\dot{\mathbf{\omega }}^{\text{b}}+\widetilde{%
\mathbf{\omega }}{^{\text{b}}}\mathbf{\Theta }^{\text{b}}\mathbf{\omega }^{%
\text{b}}+m{^{\text{b}}}\widetilde{\mathbf{d}}_{\text{bc}}\mathbf{R}^{T}%
\ddot{\mathbf{r}} &=&\mathbf{t}^{\text{b}} \\
m(\ddot{\mathbf{r}}+\mathbf{R}(\dot{\widetilde{\mathbf{\omega }}}{^{\text{b}}%
}+\widetilde{\mathbf{\omega }}^{\text{b}}\widetilde{\mathbf{\omega }}^{\text{%
b}}){^{\text{b}}}\mathbf{d}_{\text{bc}}) &=&\mathbf{f}^{\text{h}}.
\end{eqnarray}%
If a COM frame is used, combining (\ref{NEb0exp1}) and (\ref{NEhCOM2}) yields%
\begin{eqnarray}
\mathbf{\Theta }_{\text{c}}\dot{\mathbf{\omega }}_{\text{c}}^{\text{b}}+%
\widetilde{\mathbf{\omega }}_{\text{c}}^{\text{b}}\mathbf{\Theta }_{\text{c}}%
\mathbf{\omega }_{\text{c}}^{\text{b}} &=&\mathbf{t}_{\text{c}}^{\text{b}}
\label{NEcmixed} \\
m\ddot{\mathbf{r}}_{\text{c}} &=&\mathbf{f}_{\text{c}}^{\text{h}}.  \notag
\end{eqnarray}

\subsection{Arbitrary Representation}

The NE equations of body $i$ represented in an arbitrary frame $\mathcal{F}%
_{j}$ are obtained by a frame transformation of the spatial momentum balance
(\ref{Pisdot}) as%
\begin{equation}
\mathbf{Ad}_{\mathbf{C}_{j}}^{T}\dot{\mathbf{\Pi }}_{i}^{\text{s}}=\mathbf{Ad%
}_{\mathbf{C}_{j}}^{T}\mathbf{W}_{i}^{\text{s}}.
\end{equation}%
The spatial twist in terms of the twist of body $i$ represented in $\mathcal{%
F}_{j}$ is ${\mathbf{V}_{i}^{\text{s}}}=\mathbf{Ad}_{\mathbf{C}_{j}}{^{j}%
\mathbf{V}}_{i}$. Using ${\dot{\mathbf{V}}_{i}^{\text{s}}}=\mathbf{Ad}_{%
\mathbf{C}_{j}}{^{j}\dot{\mathbf{V}}}_{i}+\mathbf{ad}_{\mathbf{V}_{j}^{\text{%
s}}}\mathbf{V}_{i}^{\text{s}}$, (\ref{NEs}) yields%
\begin{equation}
{^{j}\mathbf{M}}_{i}%
\big%
({^{j}\dot{\mathbf{V}}}_{i}+\mathbf{ad}_{{^{j}\mathbf{V}}_{j}}{^{j}\mathbf{V}%
}_{i}%
\big%
)-\mathbf{ad}_{{^{j}\mathbf{V}}_{i}}^{T}{^{j}\mathbf{M}}_{i}{^{j}\mathbf{V}}%
_{i}={^{j}\mathbf{W}}_{i}  \label{NEij}
\end{equation}%
with the inertia matrix of body $i$ represented in frame $j$%
\begin{equation}
{^{j}\mathbf{M}}_{i}:=\mathbf{Ad}_{\mathbf{C}_{j}}^{T}\mathbf{M}_{i}^{\text{s%
}}\mathbf{Ad}_{\mathbf{C}_{j}}.
\end{equation}%
The spatial and body-fixed representations are special cases with $i=j$.

Even more generally, the NE equations can be resolved in yet another frame $%
\mathcal{F}_{k}$. This is achieved by transforming the momentum balance (\ref%
{Pisdot}) as%
\begin{equation}
\mathbf{Ad}_{\mathbf{R}_{j,k}}^{T}\mathbf{Ad}_{\mathbf{C}_{j}}^{T}\dot{%
\mathbf{\Pi }}_{i}^{\text{s}}=\mathbf{Ad}_{\mathbf{R}_{j,k}}^{T}\mathbf{Ad}_{%
\mathbf{C}_{j}}^{T}\mathbf{W}_{i}^{\text{s}}
\end{equation}%
where $\mathbf{R}_{k,j}$ is the rotation matrix from $\mathcal{F}_{i}$ to $%
\mathcal{F}_{k}$. The final equations follow from (\ref{NEij}) and the
relation ${^{j}\dot{\mathbf{V}}}_{i}^{j}=\mathbf{Ad}_{\mathbf{R}_{j,k}}{^{k}%
\dot{\mathbf{V}}}_{i}^{j}+\mathbf{ad}_{^{j}\overset{\omega }{\mathbf{V}}{%
_{k}^{j}}}\mathbf{Ad}_{\mathbf{R}_{j,k}}{^{k}\mathbf{V}}_{i}^{j}$ as%
\begin{equation}
\begin{tabular}{|lll|}
\hline
&  &  \\ 
& ${^{k}\mathbf{M}}_{i}^{j}%
\big%
({^{k}\dot{\mathbf{V}}}_{i}^{j}+%
\big%
(\mathbf{ad}_{{^{k}\mathbf{V}}_{j}^{j}}+\mathbf{ad}_{^{k}\overset{\omega }{%
\mathbf{V}}{_{k}^{j}}}%
\big%
){^{k}\mathbf{V}}_{i}^{j}%
\big%
)-\mathbf{ad}_{{^{k}\mathbf{V}}_{i}^{j}}^{T}{^{k}\mathbf{M}}_{i}^{j}{^{k}%
\mathbf{V}}_{i}^{j}={^{k}\mathbf{W}}_{i}^{j}$ &  \\ 
&  &  \\ \hline
\end{tabular}%
\end{equation}%
with the mass matrix of body $i$ measured at frame $\mathcal{F}_{j}$ and
resolved in frame $\mathcal{F}_{k}$%
\begin{equation}
{^{k}\mathbf{M}}_{i}^{j}:=\mathbf{Ad}_{\mathbf{R}_{k,j}}^{T}{^{j}\mathbf{M}}%
_{i}\mathbf{Ad}_{\mathbf{R}_{k,j}}=\mathbf{Ad}_{\mathbf{R}_{k,j}}^{T}\mathbf{%
Ad}_{\mathbf{C}_{j}}^{T}\mathbf{M}_{i}^{\text{s}}\mathbf{Ad}_{\mathbf{C}_{j}}%
\mathbf{Ad}_{\mathbf{R}_{k,j}}.
\end{equation}%
The spatial and body-fixed representations are special cases with $i=j=k$,
and the hybrid representation with $i=j$ and $k=0$. An alternative form of
the NE equations in arbitrary reference frames was presented in \cite%
{Bongard2015}.

\section{Recursive Evaluation of the Motion Equations for a Kinematic Chain%
\label{secRecEOM}%
}

The model-based control of complex MBS as well as the computational MBS
dynamics rely on efficient recursive inverse and forward dynamics
algorithms. A recursive Newton-Euler method for tree-topology MBS was
presented in an abstract, i.e. coordinate-free, approach in \cite{Luh1980}.
However, the various recursive methods using different representations give
rise to algorithmically equivalent methods but with different computational
costs. In the following, the various inverse dynamics algorithms are
presented and their computational effort is estimated. A detailed analysis
as well as the forward dynamics algorithms are beyond the scope of this
paper. The presented discussion is nevertheless indicative also for the
corresponding forward dynamics algorithms. Some results on the forward
kinematics complexity can be found in \cite%
{OrinSchrader1984,Stelzle1995,Yamane2009}. This depends on the actual
implementation, however. A comparative study is still due, and shall be part
of further research.

The inverse dynamics consists in evaluating the motion equations for given
joint coordinates $\mathbf{q}$, joint rates $\dot{\mathbf{q}}$,
accelerations $\ddot{\mathbf{q}}$, and applied wrenches $\mathbf{W}_{i}^{%
\text{app}}$, and hence to determine the joint forces $\mathbf{Q}=\left(
Q_{1},\ldots Q_{n}\right) $. The starting point of recursive algorithms for
rigid body MBS are the NE equations of the individual bodies. The MBS
dynamics is indeed governed by the Lagrange equations. Consequently,
summarizing the recursive steps yields the Lagrangian motion equations in
closed form. This will be shown in the following.

It is assumed for simplicity that the inertia properties, i.e. the mass
matrices $\mathbf{M}_{i}^{\text{b}}$, are expressed in the body-fixed BFR of
body $i$ determining its configuration, rather than introducing a second
frame.

\subsection{Body-fixed Representation%
\label{secInvDynBody}%
}

\paragraph{Forward Kinematics Recursion}

Given the joint variables $\mathbf{q}$, the configurations of the $n$ bodies
are determined recursively by (\ref{POEX}) or (\ref{POEY}), and the twists
by (\ref{VbRec}). Then also the accelerations are found recursively. The
expression $\mathbf{C}_{i-1,i}\left( q_{i}\right) =\mathbf{B}_{i}\exp ({^{i}}%
\mathbf{X}_{i}q_{i})$ for the relative configuration yields $\dot{\mathbf{Ad}%
}_{\mathbf{C}_{i,i-1}}\mathbf{V}_{i-1}^{\text{b}}=[\mathbf{Ad}_{\mathbf{C}%
_{i,i-1}}\mathbf{V}_{i-1}^{\text{b}},{^{i}\mathbf{X}}_{i}\dot{q}_{i}]$, and
hence%
\begin{align}
\dot{\mathbf{V}}_{i}^{\text{b}}& =\mathbf{Ad}_{\mathbf{C}_{i,i-1}}\dot{%
\mathbf{V}}_{i-1}^{\text{b}}+[\mathbf{Ad}_{\mathbf{C}_{i,i-1}}\mathbf{V}%
_{i-1}^{\text{b}},{^{i}\mathbf{X}}_{i}\dot{q}_{i}]+{^{i}\mathbf{X}}_{i}\ddot{%
q}_{i}\addtocounter{equation}{1} 
\tag{\TeXButton{\theequation}{\theequation}a}  \label{arecVbi} \\
& =\mathbf{Ad}_{\mathbf{C}_{i,i-1}}\dot{\mathbf{V}}_{i-1}^{\text{b}}+[%
\mathbf{Ad}_{\mathbf{C}_{i,i-1}}\mathbf{V}_{i-1}^{\text{b}},\mathbf{V}_{i}^{%
\text{b}}]+{^{i}\mathbf{X}}_{i}\ddot{q}_{i} 
\tag{\TeXButton{\theequation}{\theequation}b}  \label{brecVbi} \\
& =\mathbf{Ad}_{\mathbf{C}_{i,i-1}}\dot{\mathbf{V}}_{i-1}^{\text{b}}+[%
\mathbf{V}_{i}^{\text{b}},{^{i}\mathbf{X}}_{i}\dot{q}_{i}]+{^{i}\mathbf{X}}%
_{i}\ddot{q}_{i}.  \tag{\TeXButton{\theequation}{\theequation}c}
\label{crecVbi}
\end{align}%
where (\ref{brecVbi}) and (\ref{crecVbi}) follow by replacing either
argument in the Lie bracket using (\ref{VbRec}).

\begin{remark}
Notice that solving (\ref{VbRec}) for $\dot{q}_{i}$ leads to the result in
remark 9 of \cite{Part1}. Solving (\ref{crecVbi}) for $\ddot{q}_{i}$ yields (%
\ref{solqddot1}). Using (\ref{brecVbi}) the latter can be expressed as $%
\ddot{q}_{i}={^{i}\mathbf{X}}_{i}^{T}(\mathbf{V}_{i}^{\text{b}}-\mathbf{Ad}_{%
\mathbf{C}_{i,i-1}}\dot{\mathbf{V}}_{i-1}^{\text{b}}+[\mathbf{V}_{i}^{\text{b%
}},\mathbf{Ad}_{\mathbf{C}_{i,i-1}}\mathbf{V}_{i-1}^{\text{b}}])/\left\Vert {%
^{i}\mathbf{X}}_{i}\right\Vert ^{2}$.
\end{remark}

\paragraph{Recursive Newton-Euler Algorithm}

Once the configurations, twists, and accelerations of the bodies are
computed with the forward kinematics recursion, the Newton-Euler equations (%
\ref{NEb}) for each individual body can be evaluated by an inverse dynamics
backward recursion. The momentum balance of body $i$ then yields the
resulting body-fixed wrench $\mathbf{W}_{i}^{\text{b}}$ acting on the body
due to generalized joint forces and constraint reactions forces. Projecting
the resultant wrench onto the screw axis ${^{i}\mathbf{X}}_{i}$ of joint $i$
yields the generalized force $Q_{i}$. Summarizing the forward and backward
recursions yields the following recursive algorithm:%
\newpage%

\begin{itemize}
\item[ ] \underline{Forward Kinematics}%
\vspace{-0.5ex}%

\begin{itemize}
\item Input: $\mathbf{q},\dot{\mathbf{q}},\ddot{\mathbf{q}}$

\item For $i=1,\ldots ,n$ 
\vspace{-2ex}%
\begin{align}
\mathbf{C}_{i}& =\mathbf{C}_{i-1}\mathbf{B}_{i}\exp ({^{i}\mathbf{X}}%
_{i}q_{i})=\exp (\mathbf{Y}_{1}q_{1})\cdot \ldots \cdot \exp (\mathbf{Y}%
_{i}q_{i})\mathbf{A}_{i}\addtocounter{equation}{1} 
\tag{\TeXButton{\theequation}{\theequation}a}  \label{aVbRec2} \\
\mathbf{V}_{i}^{\text{b}}& =\mathbf{Ad}_{\mathbf{C}_{i,i-1}}\mathbf{V}%
_{i-1}^{\text{b}}+{^{i}\mathbf{X}}_{i}\dot{q}_{i} 
\tag{\TeXButton{\theequation}{\theequation}b}  \label{bVbRec2} \\
\dot{\mathbf{V}}_{i}^{\text{b}}& =\mathbf{Ad}_{\mathbf{C}_{i,i-1}}\dot{%
\mathbf{V}}_{i-1}^{\text{b}}-\dot{q}_{i}\mathbf{ad}_{{^{i}\mathbf{X}}_{i}}%
\mathbf{V}_{i}^{\text{b}}+{^{i}\mathbf{X}}_{i}\ddot{q}_{i} 
\tag{\TeXButton{\theequation}{\theequation}c}  \label{cVbRec2}
\end{align}%
\vspace{-4ex}%

\item Output: $\mathbf{C}_{i},\mathbf{V}_{i}^{\mathrm{b}},\dot{\mathbf{V}}%
_{i}^{\mathrm{b}}$
\end{itemize}

\item[ ] \underline{Inverse Dynamics}%
\vspace{-0.5ex}%

\begin{itemize}
\item Input: $\mathbf{C}_{i},\mathbf{V}_{i}^{\mathrm{b}},\dot{\mathbf{V}}%
_{i}^{\mathrm{b}},{\mathbf{W}}_{i}^{\text{b,app}}$

\item For $i=n-1,\ldots ,1$ 
\vspace{-2ex}
\begin{align}
\mathbf{W}_{i}^{\text{b}}& =\mathbf{Ad}_{\mathbf{C}_{i+1,i}}^{T}\mathbf{W}%
_{i+1}^{\text{b}}+\mathbf{M}_{i}^{\text{b}}\dot{\mathbf{V}}_{i}^{\text{b}}-%
\mathbf{ad}_{\mathbf{V}_{i}^{\text{b}}}^{T}\mathbf{M}_{i}^{\text{b}}\mathbf{V%
}_{i}^{\text{b}}+{\mathbf{W}}_{i}^{\text{b,app}}\addtocounter{equation}{1} 
\tag{\TeXButton{\theequation}{\theequation}a}  \label{aWbRec} \\
Q_{i}& ={^{i}\mathbf{X}}_{i}^{T}\mathbf{W}_{i}^{\text{b}} 
\tag{\TeXButton{\theequation}{\theequation}b}  \label{bWbRec}
\end{align}%
\vspace{-4ex}%

\item Output: $\mathbf{Q}$
\end{itemize}
\end{itemize}

The joint reaction wrench is omitted in (\ref{aWbRec}) since this is
reciprocal to the joint screw, and does not contribute to (\ref{bWbRec}).
Notice that, with (\ref{XY}), the body-fixed ${^{i}\mathbf{X}}_{i}$ as well
as the spatial representation $\mathbf{Y}_{i}$ of joint screw coordinates
can be used. This form of the recursive body-fixed NE equations, using Lie
group notation, has been reported in several publications \cite%
{Park1994,PloenPark1997,ParkKim2000,MUBOLieGroup}.

\paragraph{Computational Effort}

For the kinematic chain comprising $n$ bodies connected by $n$ 1-DOF joints,
in total, the twist recursion (\ref{bVbRec2}) and acceleration recursion (%
\ref{cVbRec2}) each requires $n-1$ frame transformations. The acceleration
recursion (\ref{cVbRec2}) further requires $n-1$ Lie brackets. The second
argument of the Lie bracket can be reused from (\ref{bVbRec2}). Hence the
twist and acceleration recursion need $2\left( n-1\right) $ frame
transformations and $n-1$ Lie brackets. The backward recursion (\ref{aWbRec}%
) needs $n-1$ frame transformations and $n$ Lie brackets. In total, the NE
algorithm needs $3\left( n-1\right) $ frame transformations and $2n-1$ Lie
brackets. The evaluation of the Lie bracket in (\ref{cVbRec2}) can be
simplified using (\ref{brecVbi}) since the screw vector ${^{i}\mathbf{X}}%
_{i} $ expressed in RFR is sparse and often only contains one non-zero entry.

\paragraph{Remark on Forward Dynamics}

Using the body-fixed representation, a recursive forward dynamics algorithm,
making explicit use of Lie group concepts, was presented in \cite%
{Park1994,ParkBobrowPloen1995,PloenPark1997,ParkKim2000,SohlBobrow2001}. The
kinematic forward recursion together with the factorization in section 3.1.3
of \cite{Part1} was used derive $O\left( n\right) $ forward dynamics
algorithms in \cite{Fijany1995,LillyOrin1991}, where the Lie group concept
is regarded as spatial operator algebra. Other $O\left( n\right) $ forward
dynamics algorithms were presented in \cite%
{AndersonCritchley2003,Bae2001,Hollerbach1980,HsuAnderson2002}. The inverse
dynamics formulation was also presented in \cite{Featherstone2008,Liu1988}
in the context of screw theory.

\subsection{Spatial Representation}

\paragraph{Forward Kinematics Recursion}

Expressing the spatial twist in terms of the spatial Jacobian, the
expressions (\ref{JsX}) lead immediately to%
\begin{equation}
\mathbf{V}_{i}^{\text{s}}=\mathbf{V}_{i-1}^{\text{s}}+\mathbf{J}{_{i}^{\text{%
s}}}\dot{q}_{i}.  \label{Vsirec}
\end{equation}%
The recursive determination of spatial accelerations thus only requires the
time derivative (\ref{timeDerJs}) of the spatial Jacobian, so that%
\begin{eqnarray}
\dot{\mathbf{V}}_{i}^{\text{s}} &=&\dot{\mathbf{V}}_{i-1}^{\text{s}}+\mathbf{%
ad}_{\mathbf{V}_{i}^{\text{s}}}\mathbf{J}{_{i}^{\text{s}}}\dot{q}_{i}+%
\mathbf{J}{_{i}^{\text{s}}}\ddot{q}_{i}  \label{Visdotrec} \\
&=&\dot{\mathbf{V}}_{i-1}^{\text{s}}+\mathbf{ad}_{\mathbf{V}_{i-1}^{\text{s}%
}}{\mathbf{V}{_{i}^{\text{s}}}}+\mathbf{J}{_{i}^{\text{s}}}\ddot{q}_{i}. 
\notag
\end{eqnarray}%
The second form in (\ref{Visdotrec}) follows by inserting (\ref{Vsirec}).
This is the generalization of Euler's theorem, for the derivative of vectors
resolved in moving frames, to screw coordinate vectors. Therefore the $%
\mathbf{ad}$ operator is occasionally called the 'spatial cross product'.

\paragraph{Recursive Newton-Euler Algorithm}

The momentum balance expressed with the spatial NE equations (\ref{NEs})
together with (\ref{Vsirec}) leads to the following algorithm:

\begin{itemize}
\item[ ] \underline{Forward Kinematics}%
\vspace{-0.5ex}%

\begin{itemize}
\item Input: $\mathbf{q},\dot{\mathbf{q}},\ddot{\mathbf{q}}$

\item For $i=1,\ldots ,n$ 
\vspace{-2ex}%
\begin{align}
\mathbf{C}_{i}& =\mathbf{C}_{i-1}\mathbf{B}_{i}\exp ({^{i}\mathbf{X}}%
_{i}q_{i})=\exp (\mathbf{Y}_{1}q_{1})\cdot \ldots \cdot \exp (\mathbf{Y}%
_{i}q_{i})\mathbf{A}_{i}\addtocounter{equation}{1} 
\tag{\TeXButton{\theequation}{\theequation}a}  \label{aVsRec} \\
\mathbf{J}{_{i}^{\text{s}}}& =\mathbf{Ad}_{\mathbf{C}_{i}}{^{i}\mathbf{X}}%
_{i}=\mathbf{Ad}_{\mathbf{C}_{i}\mathbf{A}_{i}^{-1}}\mathbf{Y}_{i}=\mathbf{Ad%
}_{\mathbf{C}_{j}\mathbf{S}_{j,j}}{^{j-1}}\mathbf{Z}_{j} 
\tag{\TeXButton{\theequation}{\theequation}b}  \label{bVsRec} \\
\mathbf{V}_{i}^{\text{s}}& =\mathbf{V}_{i-1}^{\text{s}}+\mathbf{J}{_{i}^{%
\text{s}}}\dot{q}_{i}  \tag{\TeXButton{\theequation}{\theequation}c}
\label{cVsRec} \\
\dot{\mathbf{V}}_{i}^{\text{s}}& =\dot{\mathbf{V}}_{i-1}^{\text{s}}+\mathbf{J%
}{_{i}^{\text{s}}}\ddot{q}_{i}+\mathbf{ad}_{\mathbf{V}_{i-1}^{\text{s}}}{%
\mathbf{V}{_{i}^{\text{s}}}}  \tag{\TeXButton{\theequation}{\theequation}d}
\label{dVsRec}
\end{align}%
\vspace{-4ex}%

\item Output: $\mathbf{C}_{i},\mathbf{V}_{i}^{\mathrm{s}},\dot{\mathbf{V}}%
_{i}^{\mathrm{s}},\mathbf{J}{_{i}^{\text{s}}}$
\end{itemize}

\item[ ] \underline{Inverse Dynamics}%
\vspace{-0.5ex}%

\begin{itemize}
\item Input: $\mathbf{C}_{i},\mathbf{V}_{i}^{\mathrm{s}},\dot{\mathbf{V}}%
_{i}^{\mathrm{s}},\mathbf{J}{_{i}^{\text{s}},\mathbf{W}}_{i}^{\text{s,app}}$

\item For $i=n-1,\ldots ,1$ 
\vspace{-2ex}
\begin{align}
\mathbf{M}_{i}^{\text{s}}& =\mathbf{Ad}_{\mathbf{C}_{i}}^{-T}\mathbf{M}_{i}^{%
\text{b}}\mathbf{Ad}_{\mathbf{C}_{i}}^{-1}\addtocounter{equation}{1} 
\tag{\TeXButton{\theequation}{\theequation}a}  \label{aWsRec} \\
\mathbf{W}_{i}^{\text{s}}& =\mathbf{W}_{i+1}^{\text{s}}+\mathbf{M}_{i}^{%
\text{s}}\dot{\mathbf{V}}_{i}^{\text{s}}-\mathbf{ad}_{\mathbf{V}_{i}^{\text{s%
}}}^{T}\mathbf{M}_{i}^{\text{s}}\mathbf{V}_{i}^{\text{s}}+{\mathbf{W}}_{i}^{%
\text{s,app}}  \tag{\TeXButton{\theequation}{\theequation}b}  \label{bWsRec}
\\
Q_{i}& =(\mathbf{J}{_{i}^{^{\text{s}}})}^{T}\mathbf{W}_{i}^{\text{s}} 
\tag{\TeXButton{\theequation}{\theequation}c}  \label{cWsRec}
\end{align}%
\vspace{-4ex}%

\item Output: $\mathbf{Q}$
\end{itemize}
\end{itemize}

\paragraph{Computational Effort}

In contrast to (\ref{bVbRec2}), once the instantaneous screws (\ref{bVsRec})
and the spatial mass matrix (\ref{Ms}) are computed, the recursions (\ref%
{cVsRec}), (\ref{dVsRec}), and (\ref{bWsRec}) do not require frame
transformations of twists. Instead the spatial mass matrix is transformed,
according to (\ref{aWsRec}), which is the frame transformation of a
second-order tensor. Overall the spatial algorithm needs $n$ frame
transformations of screw coordinates, $n$ frame transformation of a
second-order tensor, and $2n-1$ Lie brackets. Comparing body-fixed and
spatial formulation, it must be noticed that the frame transformation of the
second-order inertia tensor has the same complexity as two screw coordinate
transformations (if just implemented in the form (\ref{Ms})), and hence the
computational complexity of both would be equivalent. This fact is to be
expected since body-fixed and spatial representations are related by frame
transformations. Nevertheless the spatial version has some interesting
features that shall be emphasized:

\begin{enumerate}
\item The NE equations (\ref{NEs}) form a non-linear first-order ODE system
on $SE\left( 3\right) \times se\left( 3\right) $. Since a spatial reference
is used, the momentum conservation of a rigid body can simply be written as $%
\dot{\mathbf{\Pi }}_{i}^{\text{s}}={\mathbf{0}}$, where $\mathbf{\Pi }_{i}^{%
\text{s}}\in se^{\ast }\left( 3\right) $ is the momentum co-screw. Using the
spatial momentum balance (\ref{Pisdot}) has potentially two advantages.
Firstly, (\ref{Pisdot}) is a linear ODE in $\mathbf{\Pi }$ on the phase
space $SE\left( 3\right) \times se^{\ast }\left( 3\right) $. This implies
that a numerical integration scheme can easily preserve the momentum, as
pointed out in \cite{Bottasso1998}. Secondly, $O\left( n\right) $
formulations using canonical momenta have been shown to be computationally
advantageous. An $O\left( n\right) $ forward dynamics algorithm based on the
canonical Hamilton equations was presented in \cite{Naudet2003} that uses on
the hybrid form. It was shown to require less numerical operations than $%
O\left( n\right) $ algorithms based on the NE equations. It is also known
that $O\left( n\right) $ algorithms based on the spatial representation can
be computationally more efficient than those based on body-fixed or hybrid
representations \cite{Featherstone2008}. A further reduction of
computational costs shall be expected from an algorithm using spatial
momenta.

\item It is interesting to notice that the hybrid as well as the spatial
twists appear in the recursive $O\left( n\right) $ forward dynamics
formulation in \cite{BaeHwangHaug1988}, where the first is called 'Cartesian
velocity' and the latter 'velocity state'. In this formulation the spatial
twist plays a central role, and it was already remarked that the recursive
relation of spatial twists, ref. (\ref{cVsRec}), is simpler than that for
hybrid twists (\ref{cVhRec2}) below.\newline

\item If a \emph{purely kinematic analysis} is envisaged the forward
recursion (\ref{bVsRec})-(\ref{dVsRec}) is more efficient than the
body-fixed and the hybrid version (see next section) \cite{OrinSchrader1984}
(disregarding possibly necessary transformations of the results to local
reference frames). As pointed out in section \ref{secSpatialAcc} this
advantage is retained for the higher-order kinematics (jerk, jounce, etc.) 
\cite{MMTHighDer}.
\end{enumerate}

\paragraph{Remark on Forward Dynamics}

The spatial formulation is rarely used for dynamics. Featherstone \cite%
{Featherstone1983,Featherstone2008} derived a forward dynamics $O\left(
n\right) $ algorithm. It was concluded that this requires the lowest
computational effort compared to other methods. But this does not take into
account the necessary transformations of twists and wrenches to local
reference frames. Moreover, it was shown in \cite{Stelzle1995} that the $%
O\left( n\right) $ forward dynamics algorithm in body-fixed representation,
using the body-fixed joint screw coordinates $^{i}\mathbf{X}_{i}$ and RFR at
the joint axis, can be implemented in such a way that it requires less
computational effort than the spatial version. The key is that when the BFR
of $\mathcal{F}_{i}$ is located at and aligned with the axis of joint $i$,
then $^{i}\mathbf{X}_{i}$ becomes sparse. From a users perspective this is a
restraining presumption, however.

\subsection{Hybrid Form}

\paragraph{Forward Kinematics Recursion}

The hybrid twists is determined recursively by (\ref{VhRec}) with ${^{0}%
\mathbf{X}}_{i}^{i}{=}\mathbf{Ad}_{\mathbf{R}_{i}}{^{i}\mathbf{X}}_{i}$. For
the acceleration recursion note that $\dot{\mathbf{Ad}}_{\mathbf{r}_{i,i-1}}=%
\mathbf{ad}_{\dot{\mathbf{r}}_{i,i-1}}=\mathbf{ad}_{\dot{\mathbf{r}}_{i-1}}-%
\mathbf{ad}_{\dot{\mathbf{r}}_{i}}$ since $\dot{\mathbf{r}}_{i,i-1}=\dot{%
\mathbf{r}}_{i-1}-\dot{\mathbf{r}}_{i}$. This yields%
\begin{equation}
\dot{\mathbf{V}}_{i}^{\text{h}}=\mathbf{Ad}_{\mathbf{r}_{i,i-1}}\dot{\mathbf{%
V}}_{i-1}^{\text{h}}+{^{0}\mathbf{X}}_{i}^{i}\ddot{q}_{i}+\mathbf{ad}_{\dot{%
\mathbf{r}}_{i,i-1}}\mathbf{V}_{i-1}^{\text{h}}+\mathbf{ad}_{\mathbf{\omega }%
_{i}}{^{0}\mathbf{X}}_{i}^{i}\dot{q}_{i}.  \label{recVhdot1}
\end{equation}%
Taking into account that $\mathbf{ad}_{\dot{\mathbf{r}}_{i}}\left( \mathbf{V}%
_{i}^{\text{h}}-\mathbf{X}{_{i}^{\text{h}}}\dot{q}_{i}\right) =\mathbf{ad}_{%
\dot{\mathbf{r}}_{i}}\mathbf{V}_{i-1}^{\text{h}}$ (because there is no
angular part in $\mathbf{ad}_{\dot{\mathbf{r}}_{i}}$), and $\mathbf{ad}_{%
\dot{\mathbf{r}}_{i}}+\mathbf{ad}_{\mathbf{\omega }_{i}^{\text{s}}}=\mathbf{%
ad}_{\mathbf{V}_{i}^{\text{h}}}$, this can be transformed to%
\begin{equation}
\dot{\mathbf{V}}_{i}^{\text{h}}=\mathbf{Ad}_{\mathbf{r}_{i,i-1}}\dot{\mathbf{%
V}}_{i-1}^{\text{h}}+{^{0}\mathbf{X}}_{i}^{i}\ddot{q}_{i}+\mathbf{ad}_{\dot{%
\mathbf{r}}_{i-1}}\mathbf{V}_{i-1}^{\text{h}}-\mathbf{ad}_{\dot{\mathbf{r}}%
_{i}}\mathbf{V}_{i}^{\text{h}}+\mathbf{ad}_{\mathbf{V}_{i}^{\text{h}}}{^{0}%
\mathbf{X}}_{i}^{i}\dot{q}_{i}.  \label{recVhdot2}
\end{equation}%
Another form follows by solving (\ref{VhRec}) for $\mathbf{V}_{i-1}^{\text{h}%
}$ and inserting this into (\ref{recVhdot1}), while noting that $\mathbf{ad}%
_{\dot{\mathbf{r}}_{i,i-1}}\mathbf{Ad}_{\mathbf{r}_{i,i-1}}^{-1}=\mathbf{ad}%
_{\dot{\mathbf{r}}_{i,i-1}}$, as%
\begin{equation}
\dot{\mathbf{V}}_{i}^{\text{h}}=\mathbf{Ad}_{\mathbf{r}_{i,i-1}}\dot{\mathbf{%
V}}_{i-1}^{\text{h}}+{^{0}\mathbf{X}}_{i}^{i}\ddot{q}_{i}+\mathbf{ad}_{\dot{%
\mathbf{r}}_{i,i-1}}\mathbf{V}_{i}^{\text{h}}+(\mathbf{ad}_{\mathbf{V}_{i}^{%
\text{h}}}-\mathbf{ad}_{\dot{\mathbf{r}}_{i-1}}){^{0}\mathbf{X}}_{i}^{i}\dot{%
q}_{i}.  \label{recVhdot3}
\end{equation}%
Comparing these three different recursive relations (\ref{recVhdot1}), (\ref%
{recVhdot2}), and (\ref{recVhdot3}) for the hybrid acceleration from a
computational perspective (\ref{recVhdot1}) is the most efficient.

\paragraph{Recursive Newton-Euler Algorithm}

With the hybrid Newton-Euler equations (\ref{NEh}) the recursive NE
algorithm is as follows:

\begin{itemize}
\item[ ] \underline{Forward Kinematics}%
\vspace{-0.5ex}%

\begin{itemize}
\item Input: $\mathbf{q},\dot{\mathbf{q}},\ddot{\mathbf{q}}$

\item For $i=1,\ldots ,n$ 
\vspace{-2ex}%
\begin{align}
\mathbf{C}_{i}& =\mathbf{C}_{i-1}(\mathbf{q})\mathbf{B}_{i}\exp ({^{i}%
\mathbf{X}}_{i}q_{i})=\exp (\mathbf{Y}_{1}q_{1})\cdot \ldots \cdot \exp (%
\mathbf{Y}_{i}q_{i})\mathbf{A}_{i}\addtocounter{equation}{1} 
\tag{\TeXButton{\theequation}{\theequation}a}  \label{aVhRec2} \\
{^{0}\mathbf{X}}_{i}^{i}& =\mathbf{Ad}_{\mathbf{R}_{i}}{^{i}\mathbf{X}}_{i}=%
\mathbf{Ad}_{-\mathbf{r}_{i}}\mathbf{Y}{_{i}}=\mathbf{Ad}_{\mathbf{R}_{i}}%
\mathbf{Ad}_{\mathbf{S}_{i,i}}{^{i-1}}\mathbf{Z}_{i} 
\tag{\TeXButton{\theequation}{\theequation}b}  \label{bVhRec2} \\
\mathbf{V}_{i}^{\text{h}}& =\mathbf{Ad}_{\mathbf{r}_{i,i-1}}\mathbf{V}%
_{i-1}^{\text{h}}+{^{0}\mathbf{X}}_{i}^{i}\dot{q}_{i} 
\tag{\TeXButton{\theequation}{\theequation}c}  \label{cVhRec2} \\
\dot{\mathbf{V}}_{i}^{\text{h}}& =\mathbf{Ad}_{\mathbf{r}_{i,i-1}}\dot{%
\mathbf{V}}_{i-1}^{\text{h}}+\mathbf{ad}_{\dot{\mathbf{r}}_{i,i-1}}\mathbf{V}%
_{i-1}^{\text{h}}+\mathbf{ad}_{\mathbf{\omega }_{i}^{\text{s}}}{^{0}\mathbf{X%
}}_{i}^{i}\dot{q}_{i}+{^{0}\mathbf{X}}_{i}^{i}\ddot{q}_{i} 
\tag{\TeXButton{\theequation}{\theequation}d}  \label{dVhRec2}
\end{align}%
\vspace{-4ex}%

\item Output: $\mathbf{C}_{i},\mathbf{V}_{i}^{\mathrm{h}},\dot{\mathbf{V}}%
_{i}^{\mathrm{h}},{^{0}\mathbf{X}}_{i}^{i}$
\end{itemize}

\item[ ] \underline{Inverse Dynamics}%
\vspace{-0.5ex}%

\begin{itemize}
\item Input: $\mathbf{C}_{i},\mathbf{V}_{i}^{\mathrm{h}},\dot{\mathbf{V}}%
_{i}^{\mathrm{h}},{^{0}\mathbf{X}}_{i}^{i},{\mathbf{W}}_{i}^{\text{h,app}}$

\item For $i=n-1,\ldots ,1$ 
\vspace{-2ex}
\begin{align}
\mathbf{M}_{i}^{\text{h}}& =\mathbf{Ad}_{\mathbf{R}_{i}}^{-T}\mathbf{M}_{i}^{%
\text{b}}\mathbf{Ad}_{\mathbf{R}_{i}}^{-1}\addtocounter{equation}{1} 
\tag{\TeXButton{\theequation}{\theequation}a}  \label{aWhRec} \\
\mathbf{W}_{i}^{\text{h}}& =\mathbf{Ad}_{\mathbf{d}_{i+1,i}}^{T}\mathbf{W}%
_{i+1}^{\text{h}}+\mathbf{M}_{i}^{\text{h}}\dot{\mathbf{V}}_{i}^{\text{h}}+%
\mathbf{ad}_{\mathbf{\omega }^{\text{h}}}\mathbf{M}_{i}^{\text{h}}{\overset{%
\omega }{\mathbf{V}}}{_{i}^{\text{h}}}+{\mathbf{W}}_{i}^{\text{h,app}} 
\tag{\TeXButton{\theequation}{\theequation}b}  \label{bWhRec} \\
\mathbf{Q}_{i}& =({^{0}\mathbf{X}}_{i}^{i})^{T}\mathbf{W}_{i}^{\text{h}} 
\tag{\TeXButton{\theequation}{\theequation}c}  \label{cWhRec}
\end{align}%
\vspace{-4ex}%

\item Output: $\mathbf{Q}$
\end{itemize}
\end{itemize}

\paragraph{Computational Effort}

The hybrid representation is a compromise between using twists and wrenches
measured in body-fixed frames (as for the body-fixed representation, where
twists and wrenches are measured at the RFR origin) and those resolved in
the IFR (as for the spatial representation, where twists and wrenches are
measured at the IFR origin). It has therefore been used extensively for $%
O\left( n\right) $ inverse and forward dynamics algorithms. The essential
difference between the forward recursion for kinematic evaluation in
body-fixed and hybrid formulation is that the body-fixed recursion (\ref%
{aVbRec2}-\ref{cVbRec2}) requires frame transformations of screws involving
rotations and translations whereas the hybrid recursion (\ref{aVhRec2}-\ref%
{dVhRec2}) only requires the change of reference point using position
vectors resolved in the IFR. The attitude transformation only appears in (%
\ref{bVhRec2}) and in the computation of the hybrid inertia matrix (\ref%
{aWhRec}). In total the forward kinematics needs $n$ rotational
transformations, $2n-2$ translational transformations. Further, (\ref%
{dVhRec2}) needs $n-1$ cross products of the form $\mathbf{ad}_{\dot{\mathbf{%
r}}_{i,i-1}}\mathbf{V}_{i-1}^{\text{h}}=\left( \mathbf{0},\left( \dot{%
\mathbf{r}}_{i-1}-\dot{\mathbf{r}}_{i}\right) \times \mathbf{\omega }_{i-1}^{%
\text{s}}\right) ^{T}$ and $n$ Lie brackets $\mathbf{ad}_{\mathbf{\omega }%
_{i}^{\text{s}}}{^{0}\mathbf{X}}_{i}^{i}$. The inverse dynamics needs the $n$
rotational transformations (\ref{aWhRec}) of the second-order inertia
tensor, $n-1$ translational transformations of wrenches and $n$ Lie brackets
with $\mathbf{\omega }_{i}^{\text{s}}$ in (\ref{bWhRec}). In total the
hybrid NE algorithm needs $3n-3$ translational and $n$ rotational
transformations of screw coordinates, $n$ rotational transformations of the
inertia tensor, and $3n-1$ Lie brackets. Although the number of operations
is equivalent to the body-fixed version the particular form of
transformations is computationally very simple motivating its extensive us
in $O\left( n\right) $ forward dynamics algorithms. Moreover, the hybrid NE
equations are commonly expressed in a body-fixed BFR at the COM, so that the
hybrid NE equations simplify to (\ref{NEh0}) and (\ref{NEhCOM}),(\ref%
{NEhCOM2}), respectively.

Instead of transforming the joint screws ${^{i}\mathbf{X}}_{i}$ or $\mathbf{Y%
}{_{i}}$ in the reference configuration, the instantaneous hybrid joint
screws can be determined using the defining expression (36) in \cite{Part1}
with the current $\mathbf{b}_{j,i}$ and $\mathbf{e}_{j}$.

\paragraph{Remark on Forward Dynamics}

The above inverse dynamics formulation was presented in \cite%
{Jian1991,Jian1995,Rodriguez1991,Rodriguez1992} together with $O\left(
n\right) $ forward dynamics algorithms. An $O\left( n\right) $ forward
dynamics method was presented in \cite{Anderson1992,BaeHwangHaug1988}. These
algorithms are deemed efficient taking into account that the computation
results do not have to be transformed to the body-fixed reference points of
interest, as in case of the spatial version. An $O\left( n\right) $ forward
dynamics algorithm was developed in \cite{Naudet2003} using canonical
Hamilton equations in hybrid representation, i.e. the momentum balance (\ref%
{Pihdot}) in terms of the conjugate momenta $\mathbf{\Pi }_{i}^{\text{h}}$,
rather than the NE equations. It was concluded that its performance is
comparable to that of Featherstone's \cite{Featherstone2008} method in terms
of spatial twists.

\subsection{Choice of Body-Fixed Reference Frames%
\label{secFrameChoice}%
}

The Lie group formulation involves geometric and inertia properties that are
readily available, e.g. from CAD data.

In \cite{Part1} and in the preceding sections two approaches to the
description of MBS geometry (with and without body-fixed joint frames) and
three versions for representing velocities and accelerations (body-fixed,
spatial, hybrid) were presented, of which each has its merits. The
description of the geometry is independent from the representation of
twists. For instance, the geometry could be described in terms of joint
screws expressed in the IFR while the kinematics and dynamics is modeled
using body-fixed twists. This allows to take advantage of the low-complexity
hybrid or spatial recursive NE equations while still having the freedom to
use or avoid body-fixed joint frames.

The standard approach to model an MBS is to introduce 1.) an IFR, 2.)
body-fixed BFRs, and 3.) body-fixed JFRs. The latter is avoided using
spatial joint screws $\mathbf{Y}_{i}$, as already presented. It still
remains to introduce body-fixed BFRs kinematically representing the bodies.
However, even the \emph{explicit} definition of RFRs can be avoided by
properly placing them. Their location is usually dictated by the definition
of the inertia tensors, and it is customary to relate the inertia data to
the COM. If instead the body-fixed BFRs are assigned such that they coincide
in the reference configuration ($\mathbf{q}=\mathbf{0}$) with the IFR, then
no reference configurations of bodies need to be determined ($\mathbf{A}_{i}=%
\mathbf{I}$). This normally means that the RFR is outside the physical
extension of a body. That is, the inertia properties of all bodies are
determined in the assembly reference configuration w.r.t. to the global IFR.
In other words they are deduces from the design drawing (corresponding to $%
\mathbf{q}=\mathbf{0}$) relative to a single construction frame. This can be
exploited when using CAD systems. The required kinematic data then reduces
to the direction and position vectors, $\mathbf{e}_{i}$ and $\mathbf{y}_{i}$%
, in order to compute $\mathbf{Y}_{i}$ in (\ref{XY}). As a side effect it is 
$\mathbf{Y}_{i}={^{i}\mathbf{X}}_{i}$. This is an important result, and does
apply to any of the discussed twist representations, since the
representation of twists has nothing to do with the geometry description.
Moreover, then the POE (\ref{POEY}) and the Jacobian (\ref{JbX}), (\ref{JsX}%
), and thus (\ref{JhJb}), in terms of spatial screw coordinates simplify.
The only computational drawback is that the hybrid Newton and Euler
equations are not decoupled since the spatial IFR, to which the inertia data
is related, is unlikely to coincide with the COM of the bodies in the
reference configuration. Details can be found in \cite{IDETC2014}.

\section{Motion Equations in Closed Form%
\label{secClosedEOM}%
}

\subsection{Euler-Jourdain Equations%
\label{secProjEqu}%
}

The body-fixed NE equations for the individual bodies within the MBS are%
\begin{equation}
\mathbf{M}_{i}^{\text{b}}\dot{\mathbf{V}}_{i}^{\text{b}}-\mathbf{ad}_{%
\mathbf{V}_{i}^{\text{b}}}^{T}\mathbf{M}_{i}^{\text{b}}\mathbf{V}_{i}^{\text{%
b}}-{\mathbf{W}}_{i}^{\text{b,app}}-{\mathbf{W}}_{i}^{\text{b,c}}=\mathbf{0}
\end{equation}%
where ${\mathbf{W}}_{i}^{\text{b,c}}$ is the constraint reaction wrench of
joint $i$, and ${\mathbf{W}}_{i}^{\text{b,app}}$ represents the total wrench
applied to body $i$ including the applied wrench in joint $i$. Jourdain's
principle of virtual power, using the admissible variation $\delta {\mathsf{V%
}}^{\text{b}}=\mathsf{J}^{\text{b}}\delta \dot{\mathbf{q}}$ of the system
twist, and noting that $\delta \mathbf{V}^{\text{b}}$ are reciprocal to the
constraint wrenches (see appendix B.2), yields the system of $n$ motion
equations%
\begin{equation}
\left( \mathsf{J}^{\text{b}}\right) ^{T}\left( 
\begin{array}{c}
\mathbf{M}_{1}^{\text{b}}\dot{\mathbf{V}}_{1}^{\text{b}}-\mathbf{ad}_{%
\mathbf{V}_{1}^{\text{b}}}^{T}\mathbf{M}_{1}^{\text{b}}\mathbf{V}_{1}^{\text{%
b}}-{\mathbf{W}}_{1}^{\text{b,app}} \\ 
\vdots  \\ 
\mathbf{M}_{n}^{\text{b}}\dot{\mathbf{V}}_{n}^{\text{b}}-\mathbf{ad}_{%
\mathbf{V}_{n}^{\text{b}}}^{T}\mathbf{M}_{n}^{\text{b}}\mathbf{V}_{n}^{\text{%
b}}-{\mathbf{W}}_{n}^{\text{b,app}}%
\end{array}%
\right) =\mathbf{0}.  \label{PE}
\end{equation}%
This form allows for a concise and computationally efficient construction of
the motion equations. The point of departure are the NE equations of the
individual bodies. The body-fixed system Jacobian (\ref{JbSys}) is
determined by the (constant) joint screw coordinates in $\mathsf{X}^{\text{b}%
}$ and the screw transformations encoded in $\mathsf{A}^{\text{b}}$. The
same applies to the other representations. The accelerations are determined
by (\ref{Vbdot}), respectively (\ref{VbdotMat}). Explicit evaluation of (\ref%
{PE}) leads to the recursive algorithm in section \ref{secInvDynBody}.
Inserting the twists and accelerations in (\ref{PE}) yields the equations (%
\ref{EOM}) that determine the MBS dynamics on the tangent bundle $T{\mathbb{V%
}}^{n}$ with state vector $\left( \mathbf{q},\dot{\mathbf{q}}\right) \in T{%
\mathbb{V}}^{n}$. Alternatively, combining (\ref{PE}) with (\ref{VbAX})
yields a system of $n+6n$ ODEs in the state variables $\left( \mathbf{q},%
\mathsf{V}^{\text{b}}\right) \in {\mathbb{V}}^{n}\times se\left( 3\right)
^{n}$ that govern the dynamics on the state space ${\mathbb{V}}^{n}\times
se\left( 3\right) ^{n}$. The advantage of this formulation is that it is an
ODE system of first order, and that the system has block triangular
structure. Yet another interesting formulation follows with the NE equations
(\ref{Pisdot}) in terms of the conjugate momenta in spatial representation%
\begin{eqnarray}
\left( \mathsf{J}^{\text{s}}\right) ^{T}\left( 
\begin{array}{c}
\dot{\mathbf{\Pi }}_{1}^{\text{s}}-{\mathbf{W}}_{1}^{\text{s,app}} \\ 
\vdots  \\ 
\dot{\mathbf{\Pi }}_{n}^{\text{s}}-{\mathbf{W}}_{n}^{\text{s,app}}%
\end{array}%
\right)  &=&\mathbf{0}  \label{PEs} \\
\mathbf{M}_{i}^{\text{s}}\mathbf{J}_{i}^{\text{s}}\dot{\mathbf{q}} &=&%
\mathbf{\Pi }_{i}^{\text{s}},\ \ i,\ldots ,n.  \notag
\end{eqnarray}%
This is a system of $n+6n$ first order ODEs in the phase space $\left( 
\mathbf{q},\mathbf{\Pi }^{\text{s}}\right) \in {\mathbb{V}}^{n}\times
se^{\ast }\left( 3\right) ^{n}$. The system (\ref{PEs}) can be solved for
the $\dot{\mathbf{\Pi }}_{i}^{\text{s}}$ and $\dot{\mathbf{q}}_{i}$ noting
the block triangular structure of $\mathsf{J}^{\text{s}}$ \cite{Part1}. From
a numerical point of view the momentum formulation in phase space shall
allow for momentum preserving integration schemes.

Various versions of (\ref{PE}) have been published. Using the hybrid
representation of twists, basically the same equations were reported in \cite%
{AngelesLee1988}. There the system Jacobian is called the 'natural
orthogonal complement' motivated by the fact that the columns of $\mathsf{J}%
^{\text{b}}$ are orthogonal to the vectorial representations of constraint
wrenches (although the former are screws while the latter are co-screws). In
classical vector notation they were reported in \cite%
{IMechReport1988,maisser1988} and \cite{BremerBook2010}. In \cite%
{IMechReport1988} the equations (\ref{PE}) are called Euler-Jourdain
equations. In \cite{maisser1988}, emphasizing on the recursive evaluation of
the body Jacobian the instantaneous body-fixed joint screws $\mathbf{J}_{i}^{%
\text{b}}$ are called 'kinematic basic functions' as they are the intrinsic
objects in MBS kinematics. In \cite{BremerBook2010} the equations (\ref{PE})
are called 'projection equations' since the NE equations of the individual
bodies are restricted to the feasible motion (although $\mathsf{J}^{\text{b}%
} $ is not a projector). The equations (\ref{PE}) in body-fixed
representation are equivalent to Kane's equations where $\mathbf{J}_{i}^{%
\text{b}}$ are called 'partial velocities' \cite{KaneLevinson1985}. The
instantaneous joint screw coordinates, i.e. the columns $\mathbf{J}_{i}^{%
\text{b}}$ of the geometric Jacobian, were also called 'kinematic influence
coefficients' and their partial derivatives (\ref{der1}) the 'second-order
kinematic influence coefficients' \cite{BenedictTesar1978,ThomasTesar1982}.

It should be finally remarked that due to the block triangular form of $%
\mathsf{J}^{\text{b}}$, solving (\ref{PE}), and using the inversion of $%
\mathsf{A}^{\text{b}}$ (ref. (25) in \cite{Part1}), leads immediately to an $%
O\left( n\right) $ forward dynamics algorithm. This is the common starting
point for deriving forward dynamics algorithms that applies to any twist
representation.

\subsection{Lagrange Equations%
\label{secLagrange}%
}

The MBS motion equations can be derived as the Lagrange equations in terms
of generalized coordinates. For simplicity, potential forces are omitted so
that the Lagrangian is simply the kinetic energy. Then the equations attain
the form%
\begin{equation}
\frac{d}{dt}\left( \frac{\partial T}{\partial \dot{\mathbf{q}}}\right)
^{T}-\left( \frac{\partial T}{\partial \mathbf{q}}\right) ^{T}=\mathbf{M}%
\left( \mathbf{q}\right) \ddot{\mathbf{q}}+\mathbf{C}\left( \dot{\mathbf{q}},%
\mathbf{q}\right) \dot{\mathbf{q}}=\mathbf{Q}\left( \dot{\mathbf{q}},\mathbf{%
q},t\right)  \label{LE}
\end{equation}%
with generalized mass matrix $\mathbf{M}$, and $\mathbf{C}\left( \dot{%
\mathbf{q}},\mathbf{q}\right) \dot{\mathbf{q}}$ representing Coriolis and
centrifugal forces. The vector $\mathbf{Q}$ stands for all other generalized
forces, including potential, dissipative, and applied forces. Using
body-fixed twists, the kinetic energy of body $i$ is $T_{i}=\frac{1}{2}(%
\mathbf{V}_{i}^{\text{b}})^{T}\mathbf{M}_{i}^{\text{b}}\mathbf{V}_{i}^{\text{%
b}}$. The kinetic energy of the MBS is $T\left( \dot{\mathbf{q}},\mathbf{q}%
\right) =\sum_{i}T_{i}=\frac{1}{2}(\mathsf{V}^{\text{b}})^{T}\mathsf{M}^{%
\text{b}}\mathsf{V}^{\text{b}}=\frac{1}{2}\dot{\mathbf{q}}^{T}\mathbf{M}\dot{%
\mathbf{q}}$ with the generalized mass matrix%
\begin{equation}
\begin{tabular}{|lll|}
\hline
&  &  \\ 
& $\mathbf{M}\left( \mathbf{q}\right) =(\mathsf{J}^{\text{b}})^{T}\mathsf{M}%
^{\text{b}}\mathsf{J}^{\text{b}}$ &  \\ 
&  &  \\ \hline
\end{tabular}
\label{MLE}
\end{equation}%
and $\mathsf{M}^{\text{b}}:=\mathrm{diag}\,(\mathbf{M}_{1}^{\text{b}},\ldots
,\mathbf{M}_{n}^{\text{b}})$. The conjugate momentum vector is thus $\left( 
\frac{\partial T}{\partial \dot{\mathbf{q}}}\right) ^{T}=(\mathsf{J}^{\text{b%
}})^{T}\mathsf{M}^{\text{b}}\mathsf{V}^{\text{b}}$. Its time derivative is
with (\ref{VbdotMat}) given as $\frac{d}{dt}\left( \frac{\partial T}{%
\partial \dot{\mathbf{q}}}\right) ^{T}=\mathbf{M}\left( \mathbf{q}\right) 
\ddot{\mathbf{q}}-(\mathsf{J}^{\text{b}})^{T}\left( (\mathsf{M}^{\text{b}}%
\mathsf{A}^{\text{b}}\mathsf{a}^{\text{b}})^{T}+\mathsf{M}^{\text{b}}\mathsf{%
A}^{\text{b}}\mathsf{a}^{\text{b}}\right) \mathsf{J}^{\text{b}}\dot{\mathbf{q%
}}$, and $\mathsf{a}^{\text{b}}$ defined in (\ref{ab}). From (\ref{dJbij})
follows $\left( \frac{\partial T}{\partial \mathbf{q}}\right) ^{T}=(\mathsf{M%
}^{\text{b}}\mathsf{A}^{\text{b}}\mathsf{b}^{\text{b}}\mathsf{X}^{\text{b}%
})^{T}\mathsf{J}^{\text{b}}\dot{\mathbf{q}}$, with 
\begin{equation}
\mathsf{b}^{\text{b}}\left( \mathsf{V}^{\text{b}}\right) :=\mathrm{diag}~(%
\mathbf{ad}_{\mathbf{V}_{1}^{\text{b}}},\ldots ,\mathbf{ad}_{\mathbf{V}_{n}^{%
\text{b}}}).
\end{equation}%
This admits to identify the generalized mass matrix (\ref{MLE}) and the
matrix%
\begin{eqnarray}
\mathbf{C}\left( \mathbf{q},\dot{\mathbf{q}}\right) &=&-(\mathsf{J}^{\text{b}%
})^{T}\left( (\mathsf{M}^{\text{b}}\mathsf{A}^{\text{b}}\mathsf{a}^{\text{b}%
})^{T}+\mathsf{M}^{\text{b}}\mathsf{A}^{\text{b}}\mathsf{a}^{\text{b}%
}\right) \mathsf{J}^{\text{b}}-(\mathsf{M}^{\text{b}}\mathsf{A}^{\text{b}}%
\mathsf{b}^{\text{b}}\mathsf{X}^{\text{b}})^{T}\mathsf{J}^{\text{b}}  \notag
\\
&=&-(\mathsf{a}^{\text{b}}\mathsf{J}^{\text{b}}+\mathsf{b}^{\text{b}}\mathsf{%
X}^{\text{b}})^{T}(\mathsf{A}^{\text{b}})^{T}\mathsf{M}^{\text{b}}\mathsf{J}%
^{\text{b}}-(\mathsf{J}^{\text{b}})^{T}\mathsf{M}^{\text{b}}\mathsf{A}^{%
\text{b}}\mathsf{a}^{\text{b}}\mathsf{J}^{\text{b}}.  \label{CMatTmp}
\end{eqnarray}%
The first term on the right hand side in (\ref{CMatTmp}) can be simplified
so that%
\begin{equation}
\begin{tabular}{|lll|}
\hline
&  &  \\ 
$\mathbf{C}\left( \mathbf{q},\dot{\mathbf{q}}\right) $ & $=-(\mathsf{J}^{%
\text{b}})^{T}(\mathsf{M}^{\text{b}}\mathsf{A}^{\text{b}}\mathsf{a}^{\text{b}%
}+(\mathsf{b}^{\text{b}})^{T}\mathsf{M}^{\text{b}})\mathsf{J}^{\text{b}}.$ & 
\\ 
&  &  \\ \hline
\end{tabular}
\label{CMat}
\end{equation}%
The concise expressions (\ref{MLE}) and (\ref{CMat}) allow for construction
of the Lagrange equations in closed form. Similar expressions can be derived
using the spatial and hybrid representation of twists.

For analytic investigations of the MBS dynamics it may be useful to write
the Lagrange equations in components as%
\begin{equation}
\sum_{j=1}^{n}M_{ij}\left( \mathbf{q}\right) \ddot{q}_{j}+\sum_{j,k=1}^{n}%
\Gamma _{ijk}\left( \mathbf{q}\right) \dot{q}_{j}\dot{q}_{k}=Q_{i}\left( 
\mathbf{q},\dot{\mathbf{q}},t\right)
\end{equation}%
where the Christoffel symbols of first kind are defined as $\Gamma _{ijk}=%
\frac{1}{2}\left( \frac{\partial M_{ik}}{\partial q_{j}}+\frac{\partial
M_{ij}}{\partial q_{k}}-\frac{\partial M_{jk}}{\partial q_{i}}\right)
=\Gamma _{ikj}$. The recursive relations (\ref{der1}) give rise to the
closed form expressions 
\begin{eqnarray}
\Gamma _{ijk} &=&\frac{1}{2}\sum_{l=k}^{n}\left( (\mathbf{J}_{l,k}^{\text{b}%
})^{T}\mathbf{M}_{l}\mathbf{ad}_{\mathbf{J}_{l,i}^{\text{b}}}\mathbf{J}%
_{l,j}^{\text{b}}+(\mathbf{J}_{lj}^{\text{b}})^{T}\mathbf{M}_{l}\mathbf{ad}_{%
\mathbf{J}_{l,i}^{\text{b}}}\mathbf{J}_{l,k}^{\text{b}}+(\mathbf{J}_{l,i}^{%
\text{b}})^{T}\mathbf{M}_{l}\mathbf{ad}_{\mathbf{J}_{l,s}^{\text{b}}}\mathbf{%
J}_{l,r}^{\text{b}}\right)  \label{Gamma} \\
&&\text{with }i<j\leq k\text{ or }j\leq i<k,\ r=\max \left( i,j\right)
,s=\min \left( i,j\right) .  \notag
\end{eqnarray}%
This expression for the Christoffel symbols in Lie group notation was
reported in \cite{BrockettStokesPark1993,MUBOLieGroup}, and already in \cite%
{maisser1988} in tensor notation. This expression simplifies when Binet's
inertia tensor $\mathbf{\vartheta }_{i}=\frac{1}{2}\mathrm{tr}\,\left( 
\mathbf{\Theta }_{i}\right) \mathbf{I}-\mathbf{\Theta }_{i}$ is used in the
mass matrix $\mathbf{M}_{i}^{\text{b}}$. Then (\ref{Mbc}) is replaced by $%
\mathbf{\check{M}}_{i\text{c}}^{\text{b}}=\mathrm{diag}\,\left( \mathbf{%
\vartheta }_{i},m_{i}\mathbf{I}\right) $, and (\ref{Mb}) by $\mathbf{\check{M%
}}_{i}^{\text{b}}=\mathbf{Ad}_{\mathbf{S}_{\text{bc}}}^{-T}\mathbf{\check{M}}%
_{i\text{c}}^{\text{b}}\mathbf{Ad}_{\mathbf{S}_{\text{bc}}}^{-1}$. This
leads to%
\begin{eqnarray}
\Gamma _{ijk} &=&\frac{1}{2}\sum_{l=k}^{n}(\mathbf{J}_{l,j}^{\text{b}})^{T}%
\mathbf{\check{M}}_{l}\mathbf{ad}_{{{\overset{\omega }{\mathbf{J}}}{_{l,k}^{%
\text{b}}}}}\mathbf{J}_{l,i}^{\text{b}}  \label{GammaBinet} \\
&&\text{with }i<j\leq k\text{ or }j\leq i<k.  \notag
\end{eqnarray}

The equations (\ref{Gamma}) were presented in \cite%
{MUBOLieGroup,ParkBobrowPloen1995,ParkKim2000,Park1994,PloenPark1997}, and (%
\ref{GammaBinet}) in \cite{MUBOLieGroup}. Prior to these publications the
equations (\ref{Gamma}) and (\ref{GammaBinet}) have been reported in \cite%
{maisser1988,maisser1997} using tensor notation rather than Lie group
notation. Another publication that should be mentioned is \cite{Burdick1986}
where the Lagrange equations were derived using similar algebraic operations.

The above closed forms of EOM are derived using body-fixed twists. The
potential benefit of using spatial or hybrid twists remains to be explored.

\section{Derivatives of Motion Equations%
\label{secLinearization}%
}

In various contexts the information about the sensitivity of the MBS
kinematics and dynamics is required either w.r.t. joint angles, geometric
parameters, or dynamic parameters. Whereas it is know that the EOM of a
rigid body MBS attain a form that is linear in the dynamic parameters, they
depend non-linearly on the generalized coordinates and geometry. The POE
formulation provides a means to determine sensitivity w.r.t. to kinematic
parameters.

\subsection{Sensitivity of Motion Equations}

Gradients w.r.t. generalized coordinates are required for the linearization
of the EOM (as basis for stability analysis and controller design) as well
as for optimal control of MBS. Since the second-order and higher derivatives
(\ref{der2Jb}) of the body-fixed Jacobian, (\ref{Jsder2}) of the spatial
Jacobian, and (\ref{Jhder}) of the hybrid Jacobian are given as algebraic
closed form expressions in terms of screw products, the linearized EOM can
be evaluated recursively as well as expressed in closed-form. Using the Lie
group notation this was reported in \cite{SohlBobrow2001}.

The same results were already presented in \cite{maisser1997} using tensor
notation. Comparing the two formulations reveals once more that the matrix
Lie group formulation provides a level of abstraction leading to compact
expressions. A closed form for partial derivatives of the inverse mass
matrix has been reported in \cite{PartialDerTRO}, which is required for
investigating the controllability of MBS. Using the body-fixed
representation of twists, recursive $O\left( n\right) $ where reported in 
\cite{HsuAnderson2002,AndersonHsu2004}.

\subsection{Geometric Sensitivity}

Optimizing the design of an MBS requires information about the sensitivity
w.r.t. to geometric parameters. A recursive algorithm was reported in \cite%
{HsuAnderson2002} and its parallel implementation in \cite{AndersonHsu2004}
where the partial derivatives are computed on a case by case basis. The Lie
group formulation gives rise to a general closed-form expression. To this
end, the POE formula (\ref{POEX}) is extended as follows.

The geometry of the two bodies $i$ and $i-1$ connected by joint $i$ is
encoded in the constant part $\mathbf{S}_{i,i}$ and $\mathbf{S}_{i-1,i}$ in (%
\ref{relConf}), respectively in $\mathbf{B}_{i}$ in the formulation in (\ref%
{POEX}). These are frame transformations, and can hence be parameterized in
terms of screw coordinates. If $\mathbf{B}_{i}$ depends on $\lambda \leq 6$
geometric parameters, it is expressed as $\mathbf{B}_{i}\left( \mathbf{\pi }%
_{i}\right) =\mathbf{B}_{i0}\exp (\mathbf{U}_{i1}\pi _{i1})\cdot \ldots
\cdot \exp (\mathbf{U}_{i\lambda }\pi _{i\lambda })$. The screw coordinates $%
\mathbf{U}_{i1}$ and corresponding parameters $\pi _{i1}$ account for the
considered variations from the nominal geometry, represented by $\mathbf{B}%
_{i0}\in SE\left( 3\right) $. The relative configuration due to joint $i$
and the geometric variations is thus $\mathbf{C}_{i-1,i}\left( \mathbf{q},%
\mathbf{\pi }_{i}\right) =\mathbf{B}_{i}\left( \mathbf{\pi }_{i}\right) \exp
({^{i}}\mathbf{X}_{i}q_{i})$. The key observation is that partial
derivatives of $\mathbf{B}_{i}\left( \mathbf{\pi }_{i}\right) $ are
available in closed form, as for the joint screw coordinates. Hence also the
sensitivity w.r.t. the MBS geometry can be expressed in closed form \cite%
{MMTGenMob}. This fact has been applied to robot calibration \cite%
{Chen1997,Chen2001} where the POE accounts for geometric imperfections to be
identified.

\subsection{Time Derivatives of the EOM}

The design of feedback-linearizing flatness-based controllers for robotic
manipulators that are modeled as rigid body MBS actuated by elastic
actuators (so-called series elastic actuators) requires the time derivatives
of the inverse dynamics solution $\mathbf{Q}\left( t\right) $ \cite%
{deLuca1998,PalliMelchiorriDeLuca2008}. That is, the first and second time
derivatives of the EOM are necessary. Extensions of the classical recursive
Newton-Euler inverse dynamics algorithms in body-fixed representations were
presented in \cite{BuondonnaDeLuca2015}. As it can be expected the relation
are very complicated. Using the presented Lie group formulation of the
inverse dynamics algorithms gives rise to rather compact and thus fail-safe
algorithm. This was presented in \cite{ICRA2017} for the body-fixed and
hybrid version.

\section{Geometric Integration%
\label{secIntegration}%
}

This paper focussed on the MBS modeling in terms of relative (joint)
coordinates. Alternatively, the MBS kinematics can be described in terms of
absolute coordinates.

One of the issues that is being addressed when modeling MBS in terms of
absolute coordinates is the \emph{kinematic reconstruction}, i.e. the
determination of the motion of a rigid body, represented by $\mathbf{C}%
\left( t\right) $, from its velocity field $\mathbf{V}\left( t\right) $.
This amounts to solving one of the equations (see appendix A2 in \cite{Part1}%
) 
\begin{equation}
\widehat{\mathbf{V}}^{\text{b}}=\mathbf{C}^{-1}\dot{\mathbf{C}},\ \ \ \ \ \
\ \ \widehat{\mathbf{V}}^{\text{s}}=\dot{\mathbf{C}}\mathbf{C}^{-1}
\label{Vbs}
\end{equation}%
together with the NE (\ref{NEb}) or (\ref{NEs}), respectively. Classically,
the orientation is parameterized with three parameters. The problem
encountered is that there is no singularity-free global parameterization of
rotations with three parameters. Instead of local parameters (position and
rotation angles) the absolute configurations of the rigid bodies within the
MBS can be represented by $\mathbf{C}\left( t\right) $. Then a numerical
integration step from time $t_{k-1}$ to $t_{k}=t_{k-1}+h$ shall determine
the incremental configuration update $\Delta \mathbf{C}_{k}=\mathbf{C}%
_{k-1}^{-1}\mathbf{C}_{k}$ with $\mathbf{C}_{k}=\mathbf{C}\left(
t_{k}\right) $ and $\mathbf{C}_{k-1}=\mathbf{C}\left( t_{k-1}\right) $. The
equations (\ref{Vbs}) are ODEs on the Lie group $SE\left( 3\right) $. These
can be replaced by ODEs on the Lie algebra $se\left( 3\right) $. The motion
increment from $t_{k-1}$ to $t_{k}$ is parameterized as $\Delta \mathbf{C}%
\left( t\right) =\exp \mathbf{X}\left( t\right) $ with an algorithmic
instantaneous screw coordinate vector $\mathbf{X}$. Then (\ref{Vbs}) are
equivalent to the ODEs on the Lie algebra 
\begin{equation}
\mathbf{V}^{\text{s}}=\mathbf{dexp}_{\mathbf{X}}\dot{\mathbf{X}},\ \ \ \ 
\mathbf{V}^{\text{b}}=\mathbf{dexp}_{-\mathbf{X}}\dot{\mathbf{X}}
\label{Vdexp}
\end{equation}%
where $\mathbf{dexp}_{\mathbf{X}}:se\left( 3\right) \rightarrow se\left(
3\right) $ is the right-trivialized differential of the exp mapping on $%
SE\left( 3\right) $ \cite{Bottasso1998,BIT,CND2016,ParkChung2005}. This is
the basic idea of the class of Munthe-Kaas integration schemes \cite%
{CelledoniOwren1999,Iserles2000,muntekaas1999}. This scheme has been adapted
to MBS in absolute coordinates \cite{MUBO2014}. The advantage of these
integration methods is that no global parameterization is necessary since
the numerical integration is pursued in terms of the incremental parameters $%
\mathbf{X}$. The ODEs (\ref{Vdexp}) can be solved with any vector space
integration scheme (originally the Munthe-Kaas scheme uses a Runge-Kutta
method) with initial value $\mathbf{X}\left( t_{k-1}\right) =\mathbf{0}$.

Recently the geometric integration concepts were incorporated in the
generalized $\alpha $ method \cite{Krysl,BruelsCardonaArnold2012} for MBS
described in absolute coordinates. In this case the representation of proper
rigid body motions is crucial as discussed in \cite{BIT,MMT2014}, which is
frequently incorrectly represented by $SO\left( 3\right) \times {\mathbb{R}}%
^{3}$. Also momentum preserving schemes were proposed \cite{CND2015}. It
should be mentioned that the concept of geometric integration schemes on $%
SE\left( 3\right) $ can be transferred to the kinematics of flexible bodies
undergoing large deformations described as Cosserat continua. In this
context the spatial description (referred to as fixed pole formulation) has
proven to be beneficial \cite{GacesaJelenic2015}. Recent results on Lie
group modeling of beams can be found in \cite%
{SonnevilleCardonaBruls14,SonnevilleBruls14}.

\section{Conclusions and Outlook}

The computational effort of recursive $O\left( n\right) $ algorithms, but
also of the formalisms for evaluating the EOM in closed form, depends on the
representation of rigid body motions and of the motions of technical joints.
Since the geometry of finite rigid body and relative motions is described by
the Lie group $SE\left( 3\right) $, and that of instantaneous motions be the
screw algebra $se\left( 3\right) $, Lie group theory provides the geometric
framework. As already shown in \cite{Part1}, Lie group formulations for the
MBS kinematics give rise to compact recursive formulations in terms of
relative coordinates. In this paper the corresponding recursive NE
algorithms were presented and related to the various $O\left( n\right) $
algorithms scattered in the literature. This allows for a comparative
investigation of their efficiency in conjunction with the modeling
procedure. For instance, whereas most $O\left( n\right) $ algorithms used
the hybrid representation, the spatial representation, as used by
Featherstone \cite{Featherstone2008} and Bottasso \cite{Bottasso1998} (where
it is called fixed point formulation), is receiving increased attention
since it gives easily rise to structure preserving integration schemes \cite%
{Borri2001b,Borri2003,Bottasso1998,GacesaJelenic2015}. A conclusive
investigation will be the subject of future research. Future research will
also focus on combining the $O\left( n\right) $ forward dynamics algorithm
by Featherstone \cite{Featherstone2008}, based on NE equations using spatial
representations with Naudet's algorithm \cite{Naudet2003} based on
Hamilton's canonical equations in hybrid representation. The use of the
spatial momentum balance shall allow for momentum preserving integration of
the EOM and at the same time to reduce the number of frame transformations.
A further important research topic is the derivation of structure preserving
Lie group integration schemes for which the spatial formulation of EOM will
be formulation of choice.

\section*{A Summary of basic Kinematic Relations}

As prerequisite the kinematic relations derived in \cite{Part1} are
summarized. Denote with%
\begin{equation}
\mathbf{C}_{i}=\left( 
\begin{array}{cc}
\mathbf{R}_{i} & \mathbf{r}_{i} \\ 
\mathbf{0} & 1%
\end{array}%
\right) \in SE\left( 3\right)
\end{equation}%
the \emph{absolute configuration} of body $i$ w.r.t. the inertial frame
(IFR) $\mathcal{F}_{0}$. This is alternatively denoted with $C_{i}=\left( 
\mathbf{R}_{i},\mathbf{r}_{i}\right) $. The \emph{relative configuration} of
body $i$ relative to body $i-1$ is given as%
\begin{equation}
\mathbf{C}_{i-1,i}\left( q_{i}\right) =\mathbf{S}_{i-1,i}\exp ({}{}{^{i-1}}%
\mathbf{Z}_{i}q_{i})\mathbf{S}_{i,i}^{-1}=\mathbf{B}_{i}\exp ({^{i}}\mathbf{X%
}_{i}q_{i})  \label{relConf}
\end{equation}%
where $\mathbf{B}_{i}:=\mathbf{S}_{i-1,i}\mathbf{S}_{i,i}^{-1}=\mathbf{C}%
_{i-1,i}\left( 0\right) $ is the reference configuration of body $i$ w.r.t.
body $i-1$, i.e. for $q_{i}=0$, and ${^{i-1}}\mathbf{Z}_{i}\in {\mathbb{R}}%
^{6}$ is the screw coordinate vector of joint $i$ represented in the joint
frame (JFR) $\mathcal{J}_{i-1,i}$ on body $i-1$. Successive relative
configurations can be combined to%
\begin{eqnarray}
\mathbf{C}_{i}\left( \mathbf{q}\right) &=&\mathbf{B}_{1}\exp ({^{1}}\mathbf{X%
}_{1}q_{1})\cdot \mathbf{B}_{2}\exp ({^{2}}\mathbf{X}_{2}q_{2})\cdot \ldots
\cdot \mathbf{B}_{i}\exp ({^{i}}\mathbf{X}_{i}q_{i})  \label{POEX} \\
&=&\exp (\mathbf{Y}_{1}q_{1})\cdot \exp (\mathbf{Y}_{2}q_{2})\cdot \ldots
\cdot \exp (\mathbf{Y}_{i}q_{i})\mathbf{A}_{i}  \label{POEY}
\end{eqnarray}%
where ${^{i}}\mathbf{X}_{i}\in {\mathbb{R}}^{6}$ is the screw coordinate
vector of joint $i$ represented in the joint frame fixed at body $i$, $%
\mathbf{Y}_{i}\in {\mathbb{R}}^{6}$ is the joint screw coordinate vector in
spatial representation (measured and resolved in IFR) for the reference
configuration $\mathbf{q}=\mathbf{0}$, and $\mathbf{A}_{i}=\mathbf{C}%
_{i}\left( \mathbf{0}\right) $ is the reference configuration of body $i$.
The two representations of joint screw coordinates are related by%
\begin{equation}
\mathbf{Y}_{i}=\mathbf{Ad}_{\mathbf{A}_{i}}{^{i}}\mathbf{X}_{i},\ \ {^{i}}%
\mathbf{X}_{i}=\mathbf{Ad}_{\mathbf{S}_{i,i}}{^{i-1}}\mathbf{Z}_{i}
\label{XY}
\end{equation}%
where, in vector representation of screws, the adjoined transformation $%
\mathbf{Ad}$ corresponding to $\mathbf{C}\in SE\left( 3\right) $ is given by
the matrix%
\begin{equation}
\mathbf{Ad}_{\mathbf{C}}=\left( 
\begin{array}{cc}
\mathbf{R} & \mathbf{0} \\ 
\widetilde{\mathbf{r}}\mathbf{R}\ \  & \mathbf{R}%
\end{array}%
\right) .  \label{Ad}
\end{equation}%
For sake of simplicity, the following notations are used%
\begin{equation}
\mathbf{Ad}_{\mathbf{R}}=\left( 
\begin{array}{cc}
\mathbf{R} & \mathbf{0} \\ 
\mathbf{0} & \mathbf{R}%
\end{array}%
\right) ,\ \text{for}\ \ C=\left( \mathbf{R},\mathbf{0}\right) \ \ \ \ \ \ \
\ \mathbf{Ad}_{\mathbf{r}}=\left( 
\begin{array}{cc}
\mathbf{I} & \mathbf{0} \\ 
\widetilde{\mathbf{r}} & \mathbf{I}%
\end{array}%
\right) ,\ \text{for}\ \ C=\left( \mathbf{I},\mathbf{r}\right)  \label{AdRr}
\end{equation}%
so that $\mathbf{Ad}_{\mathbf{C}}=\mathbf{Ad}_{\mathbf{r}}\mathbf{Ad}_{%
\mathbf{R}}$.

The twist of body $i$ in \emph{body-fixed} representation $\mathbf{V}_{i}^{%
\text{b}}=\left( \mathbf{\omega }_{i}^{\text{b}},\mathbf{v}_{i}^{\text{b}%
}\right) ^{T}$ and in \emph{spatial} representation $\mathbf{V}_{i}^{\text{s}%
}=\left( \mathbf{\omega }_{i}^{\text{s}},\mathbf{v}_{i}^{\text{s}}\right)
^{T}$ is defined by%
\begin{equation}
\widehat{\mathbf{V}}_{i}^{\text{b}}=\left( 
\begin{array}{cc}
\widetilde{\mathbf{\omega }}_{i}^{\text{b}} & \mathbf{v}_{i}^{\text{b}} \\ 
\mathbf{0} & 0%
\end{array}%
\right) =\mathbf{C}_{i}^{-1}\dot{\mathbf{C}}_{i},\ \widehat{\mathbf{V}}_{i}^{%
\text{s}}=\left( 
\begin{array}{cc}
\widetilde{\mathbf{\omega }}_{i}^{\text{s}} & \mathbf{v}_{i}^{\text{s}} \\ 
\mathbf{0} & 0%
\end{array}%
\right) =\dot{\mathbf{C}}_{i}\mathbf{C}_{i}^{-1}.
\end{equation}%
Here $\mathbf{v}_{i}^{\text{b}}=\mathbf{R}_{i}^{T}\dot{\mathbf{r}}_{i}$ the
body-fixed translational velocity, i.e. the velocity of the origin of the
body-fixed reference frame (RFR) $\mathcal{F}_{i}$ of body $i$ measured in
the IFR $\mathcal{F}_{0}$ and resolved in $\mathcal{F}_{i}$, whereas $%
\mathbf{v}_{i}^{\text{s}}=\dot{\mathbf{r}}_{i}+\mathbf{r}_{i}\times \mathbf{%
\omega }_{i}^{\text{s}}$ is the spatial translational velocity, i.e. the
velocity of the point of the body that is momentarily passing through the
origin of the IFR $\mathcal{F}_{0}$ resolved in the IFR. The body-fixed and
spatial angular velocity, $\mathbf{\omega }_{i}^{\text{b}}$ and $\mathbf{%
\omega }_{i}^{\text{s}}$, is defined by $\widetilde{\mathbf{\omega }}_{i}^{%
\text{b}}=\mathbf{R}_{i}^{T}\dot{\mathbf{R}}_{i}$ and $\widetilde{\mathbf{%
\omega }}_{i}^{\text{s}}=\dot{\mathbf{R}}_{i}\mathbf{R}_{i}^{T}$,
respectively. The \emph{hybrid twist} is defined as $\mathbf{V}_{i}^{\text{h}%
}=(\mathbf{\omega }_{i}^{\text{s}},\dot{\mathbf{r}}_{i})^{T}$, and finally
the \emph{mixed twist} as $\mathbf{V}_{i}^{\text{m}}=(\mathbf{\omega }_{i}^{%
\text{b}},\dot{\mathbf{r}}_{i})^{T}$. The four representations are related
as follows%
\begin{eqnarray}
\mathbf{V}_{i}^{\text{h}} &=&\left( 
\begin{array}{cc}
\mathbf{R}_{i} & \mathbf{0} \\ 
\mathbf{0} & \mathbf{R}_{i}%
\end{array}%
\right) \mathbf{V}_{i}^{\text{b}}=\mathbf{Ad}_{\mathbf{R}_{i}}\mathbf{V}%
_{i}^{\text{b}}  \label{VhVb} \\
\mathbf{V}_{i}^{\text{s}} &=&\mathbf{Ad}_{\mathbf{C}_{i}}\mathbf{V}_{i}^{%
\text{b}}=\mathbf{Ad}_{\mathbf{C}_{i}}\mathbf{Ad}_{\mathbf{R}_{i}}^{-1}%
\mathbf{V}_{i}^{\text{h}}=\mathbf{Ad}_{\mathbf{r}_{i}}\mathbf{V}_{i}^{\text{h%
}}  \label{VsVb} \\
\mathbf{V}_{i}^{\text{m}} &=&\left( 
\begin{array}{cc}
\mathbf{I} & \mathbf{0} \\ 
\mathbf{0} & \mathbf{R}_{i}%
\end{array}%
\right) \mathbf{V}_{i}^{\text{b}}=\left( 
\begin{array}{cc}
\mathbf{R}_{i}^{T} & \mathbf{0} \\ 
\mathbf{0} & \mathbf{I}%
\end{array}%
\right) \mathbf{V}_{i}^{\text{h}}=\left( 
\begin{array}{cc}
\mathbf{R}_{i}^{T} & \mathbf{0} \\ 
-\widetilde{\mathbf{r}}_{i} & \mathbf{I}%
\end{array}%
\right) \mathbf{V}_{i}^{\text{s}}.  \label{Vmbhs}
\end{eqnarray}%
The twist of body $i$ within a kinematic chain is determined in terms of the
generalized velocities $\dot{\mathbf{q}}$ as%
\begin{eqnarray}
\mathbf{V}_{i}^{\text{b}} &=&\sum_{j\leq i}\mathbf{J}_{i,j}^{\text{b}}\dot{q}%
_{j}=\mathsf{J}_{i}^{\text{b}}\dot{\mathbf{q}},\ \ \ \ \ \ \ \mathbf{V}_{i}^{%
\text{s}}=\sum_{j\leq i}\mathbf{J}_{j}^{\text{s}}\dot{q}_{j}=\mathsf{J}_{i}^{%
\text{s}}\dot{\mathbf{q}}  \label{VbJac} \\
\mathbf{V}_{i}^{\text{h}} &=&\sum_{j\leq i}\mathbf{J}_{i,j}^{\text{h}}\dot{q}%
_{j}=\mathsf{J}_{i}^{\text{h}}\dot{\mathbf{q}},\ \ \ \ \ \ \ \mathbf{V}_{i}^{%
\text{m}}=\sum_{j\leq i}\mathbf{J}_{i,j}^{\text{m}}\dot{q}_{j}=\mathsf{J}%
_{i}^{\text{m}}\dot{\mathbf{q}}
\end{eqnarray}%
with the Jacobian $\mathsf{J}_{i}^{\text{b}}$ in body-fixed, $\mathsf{J}%
_{i}^{\text{s}}$ in spatial, $\mathsf{J}_{i}^{\text{h}}$ in hybrid, and $%
\mathsf{J}_{i}^{\text{m}}$ in mixed representation. The $i$th column of the
Jacobian is respectively given by%
\begin{eqnarray}
\mathbf{J}_{i,j}^{\text{b}} &=&\mathbf{Ad}_{\mathbf{C}_{i,j}}{^{j}}\mathbf{X}%
_{j}=\mathbf{Ad}_{\mathbf{C}_{i,j}\mathbf{A}_{j}^{-1}}\mathbf{Y}_{j}
\label{JbX} \\
\mathbf{J}_{j}^{\text{s}} &=&\mathbf{Ad}_{\mathbf{C}_{j}}{^{j}}\mathbf{X}%
_{j}=\mathbf{Ad}_{\mathbf{C}_{j}\mathbf{A}_{j}^{-1}}\mathbf{Y}_{j}
\label{JsX} \\
\mathbf{J}_{i,j}^{\text{h}} &=&\mathbf{Ad}_{\mathbf{r}_{i,j}}{^{0}}\mathbf{X}%
_{j}^{j},\ \ \ \ \ \ \ \ \ \ \ \ \ \ \ \ \ \ \ \ \ \ \text{for }j\leq i.
\label{Jh}
\end{eqnarray}%
These are the instantaneous joint screw coordinates in body-fixed, spatial,
and hybrid representation. The Jacobians are related as%
\begin{equation}
\mathbf{J}_{j}^{\text{s}}=\mathbf{Ad}_{\mathbf{r}_{i}}\mathbf{J}_{i,j}^{%
\text{h}}=\mathbf{Ad}_{\mathbf{C}_{i}}\mathbf{J}_{i}^{\text{b}},\ \ \mathbf{J%
}_{i}^{\text{h}}=\mathbf{Ad}_{\mathbf{R}_{i}}\mathbf{J}_{i}^{\text{b}}.
\label{JhJb}
\end{equation}%
The representations of joint screw coordinates are related by (\ref{XY}) and
by%
\begin{equation}
\mathbf{Y}_{j}=\mathbf{Ad}_{\mathbf{r}_{i}}{^{0}}\mathbf{X}_{j}^{j},\ \ \ {%
^{0}}\mathbf{X}_{j}^{j}=\mathbf{Ad}_{\mathbf{R}_{j}}{^{j}}\mathbf{X}_{j}
\label{YsXh}
\end{equation}%
where $\mathbf{r}_{i}$ is the current position of body $i$ in $\mathbf{C}%
_{i} $. The twists admit the recursive expressions 
\begin{eqnarray}
\mathbf{V}_{i}^{\text{b}} &=&\mathbf{Ad}_{\mathbf{C}_{i,i-1}}\mathbf{V}%
_{i-1}^{\text{b}}+{^{i}}\mathbf{X}_{i}\dot{q}_{i}  \label{VbRec} \\
\mathbf{V}_{i}^{\text{s}} &=&\mathbf{V}_{i-1}^{\text{s}}+\mathbf{J}{_{i}^{%
\text{s}}}\dot{q}_{i}  \label{Vsrec} \\
\mathbf{V}_{i}^{\text{h}} &=&\mathbf{Ad}_{\mathbf{r}_{i,i-1}}\mathbf{V}%
_{i-1}^{\text{h}}+{^{0}}\mathbf{X}_{i}^{i}\dot{q}^{i}.  \label{VhRec}
\end{eqnarray}%
Summarizing the twists of all bodies in $\mathsf{V}^{\text{b}},\mathsf{V}^{%
\text{s}},\mathsf{V}^{\text{h}}\in {\mathbb{R}}^{6n}$, respectively, admits
the expressions 
\begin{equation}
\mathsf{V}^{\text{b}}=\mathsf{J}^{\text{b}}\dot{\mathbf{q}},\ \ \mathsf{V}^{%
\text{s}}=\mathsf{J}^{\text{s}}\dot{\mathbf{q}},\ \ \mathsf{V}^{\text{h}}=%
\mathsf{J}^{\text{h}}\dot{\mathbf{q}}  \label{VbAX}
\end{equation}%
in terms of the system Jacobians that admit the factorizations \cite{Part1}%
\begin{equation}
\mathsf{J}^{\text{b}}=\mathsf{A}^{\text{b}}\mathsf{X}^{\text{b}},\ \ \mathsf{%
J}^{\text{s}}=\mathsf{A}^{\text{s}}\mathsf{Y}^{\text{s}}=\mathsf{A}^{\text{sb%
}}\mathsf{X}^{\text{b}},\ \ \mathsf{J}^{\text{h}}=\mathsf{A}^{\text{h}}%
\mathsf{X}^{\text{h}}.  \label{JbSys}
\end{equation}%
This provides a compact description of the overall MBS kinematics. The
explicit relations for the inverse of the matrices $\mathsf{A}$ is the
starting point for deriving recursive forward dynamics $O\left( n\right) $
algorithms.

\section*{B Rigid Body Motions and the Lie Group $SE\left( 3\right) $}

For an introduction to screws and to the motion Lie group $SE\left( 3\right) 
$ the reader is referred to the text books \cite%
{Angeles2003,LynchPark2017,Murray,Selig}.

\subsection*{B.1 Derivatives of Screws}

Let $\mathbf{C}_{i}$ be time dependent. According to (\ref{Ad}), the
corresponding frame transformation of screw coordinates from $\mathcal{F}%
_{i} $ to $\mathcal{F}_{0}$ is ${\mathbf{X}}\equiv {^{0}\mathbf{X}}=\mathbf{%
Ad}_{\mathbf{C}_{i}}{^{i}}\mathbf{X}$. Assume that the screw coordinates
expressed in body-fixed frame are constant. The rate of change of the screw
coordinates expressed in IFR is $\frac{d}{dt}\widehat{\mathbf{X}}=\frac{d}{dt%
}\mathbf{Ad}_{\mathbf{C}_{i}}({^{i}\widehat{\mathbf{X}})}=\dot{\mathbf{C}}%
_{i}\mathbf{C}_{i}^{-1}\mathbf{C}_{i}{^{i}\widehat{\mathbf{X}}}\mathbf{C}%
_{i}^{-1}-\mathbf{C}_{i}{^{i}\widehat{\mathbf{X}}}\mathbf{C}_{i}^{-1}\dot{%
\mathbf{C}}_{i}\mathbf{C}_{i}^{-1}=\widehat{\mathbf{V}}_{i}^{\text{s}}\,{%
\widehat{\mathbf{X}}}-\widehat{\mathbf{V}}_{i}^{\text{s}}\,{\widehat{\mathbf{%
X}}}=[\widehat{\mathbf{V}}_{i}^{\text{s}},{\widehat{\mathbf{X}}}]$. Therein 
\begin{equation}
\lbrack {\widehat{\mathbf{X}}}_{1}{,\widehat{\mathbf{X}}}_{2}]={\widehat{%
\mathbf{X}}}_{1}{\widehat{\mathbf{X}}}_{2}-{\widehat{\mathbf{X}}}_{2}{%
\widehat{\mathbf{X}}}_{1}=\mathrm{ad}_{\mathbf{X}_{1}}(\mathbf{X}_{2})
\label{LieBracketSE3}
\end{equation}%
is the Lie bracket on $se\left( 3\right) $, also called the \emph{adjoint}
mapping. In vector notation of screws, denoting a general screw vector with $%
\mathbf{X}=\left( \mathbf{\xi },\mathbf{\eta }\right) ^{T}$, this is%
\begin{equation}
\lbrack {\mathbf{X}}_{1}{,\mathbf{X}}_{2}]=\left( \mathbf{\xi }_{1}\times 
\mathbf{\xi }_{2},\mathbf{\eta }_{1}\times \mathbf{\xi }_{2}+\mathbf{\xi }%
_{1}\times \mathbf{\eta }_{2}\right) ^{T}=\mathbf{ad}_{\mathbf{X}_{1}}%
\mathbf{X}_{2}  \label{ScrewProd}
\end{equation}%
with 
\begin{equation}
\mathbf{ad}_{\mathbf{X}}=\left( 
\begin{array}{cc}
\widetilde{\mathbf{\xi }}\ \  & \mathbf{0} \\ 
\widetilde{\mathbf{\eta }}\ \  & \widetilde{\mathbf{\xi }}%
\end{array}%
\right) .  \label{adse3}
\end{equation}%
The form (\ref{ScrewProd}) is known as the screw product \cite%
{BottemaRoth1990,Selig}. The matrix (\ref{adse3}) has appeared under
different names, such as 'spatial cross product' in \cite%
{Featherstone2001,Featherstone2008,Jian1995}, or the 'north-east cross
product' \cite{Bottasso1998}. The Lie bracket obeys the Jacobi identity%
\begin{equation}
\lbrack {\mathbf{X}}_{1}{,[\mathbf{X}}_{2}\mathbf{,}{\mathbf{X}}_{3}]+[{%
\mathbf{X}}_{2}{,[\mathbf{X}}_{3}\mathbf{,}{\mathbf{X}}_{1}]+[{\mathbf{X}}%
_{3}{,[\mathbf{X}}_{1}\mathbf{,}{\mathbf{X}}_{2}]=\mathbf{0}.
\label{JacIdent}
\end{equation}%
Allowing for time dependent body-fixed screw coordinates ${^{i}\mathbf{X}}$,
the above relation gives rise to an expression for the time derivative of
screw coordinates in moving frames%
\begin{equation}
\mathbf{Ad}_{\mathbf{C}_{i}}^{-1}\dot{\mathbf{X}}={^{i}}\dot{\mathbf{X}}+[%
\mathbf{V}_{i}^{\text{b}},{^{i}\mathbf{X}}].
\end{equation}%
This is the spatial extension of Euler's formula for the derivative of a
vector resolved in a moving frame.

For the sake of simplicity throughout the paper, the following notations are
used%
\begin{equation}
{{\overset{\omega }{\mathbf{V}}}}=\left( 
\begin{array}{c}
\mathbf{\omega } \\ 
\mathbf{0}%
\end{array}%
\right) ,\ \ \ {{\overset{v}{\mathbf{V}}}}=\left( 
\begin{array}{c}
\mathbf{0} \\ 
\mathbf{v}%
\end{array}%
\right) .  \label{Vw}
\end{equation}%
Then the matrices%
\begin{equation}
\mathbf{ad}_{\mathbf{\omega }}=\left( 
\begin{array}{cc}
\widetilde{\mathbf{\omega }}\ \  & \mathbf{0} \\ 
\mathbf{0}\ \  & \widetilde{\mathbf{\omega }}%
\end{array}%
\right) ,\ \ \ \mathbf{ad}_{\mathbf{v}}=\left( 
\begin{array}{cc}
\mathbf{0}\  & \mathbf{0} \\ 
\widetilde{\mathbf{v}}\ \  & \mathbf{0}%
\end{array}%
\right)
\end{equation}%
are used to denote the matrix (\ref{adse3}) for twists ${{\overset{\omega }{%
\mathbf{V}}}}$ and ${{\overset{v}{\mathbf{V}}}}$, respectively, which are
the infinitesimal versions on (\ref{AdRr}).

\subsection*{B.2 Wrenches as Co-Screws -- $se^{\ast }\left( 3\right) $}

Screws are the geometric objects embodying twists, wrenches, and momenta of
rigid bodies. These different physical meanings imply different mathematical
interpretations of the geometric object.

A wrench, defined by a force and moment, is denoted with $\mathbf{W}=\left( 
\mathbf{t},\mathbf{f}\right) ^{T}$. The force applied at a point with
position vector $\mathbf{p}$ generates the moment $\mathbf{t}=\mathbf{p}%
\times \mathbf{f}$. The dual to Chasles theorem is the Poinsot theorem
stating that every system of forces can be reduced to a force together with
a couple with moment parallel to the force.

Geometrically a screw is determined by the Pl\"{u}cker coordinates of the
line along the screw axis and the pitch. If $\mathbf{e}$ is the unit vector
along the screw axis, and $\mathbf{p}$ is a position vector of a point on
that axis, the screw coordinate vector of a twist is $\mathbf{V}=\left( 
\mathbf{\omega },\mathbf{v}\right) ^{T}=\omega \left( \mathbf{e},\mathbf{p}%
\times \mathbf{e}+h\mathbf{e}\right) ^{T}$, where $\omega =\left\Vert 
\mathbf{\omega }\right\Vert $ is its magnitude, and $h=\mathbf{v}^{T}\mathbf{%
\omega }/\omega ^{2}$ is its pitch. The screw coordinate vector of a wrench,
i.e. the force $\mathbf{f}$ producing a torque $\mathbf{t}$ about the axis $%
\mathbf{e}$ when the point of application is displaced according to $\mathbf{%
p}$ from the axis, is $\mathbf{W}=\left( \mathbf{t},\mathbf{f}\right)
^{T}=f\left( \mathbf{p}\times \mathbf{e}+h\mathbf{e},\mathbf{e}\right) ^{T}$%
, with pitch $h=\mathbf{t}^{T}\mathbf{f}/\left\Vert \mathbf{f}\right\Vert
^{2}$. Apparently the linear and angular components of the screw coordinates
are interchanged for twists and wrenches. The different definition of screw
coordinate vectors allows to describe the action of a wrench on a twist as
the scalar product: $\mathbf{W}^{T}\mathbf{V}$ is the power performed by the
wrench acting on twist $\mathbf{V}$.

A twist ${^{2}}\mathbf{V}$ represented in frame $\mathcal{F}_{2}$ transforms
to its representation in frame $\mathcal{F}_{1}$ according to ${^{1}}\mathbf{%
V}=\mathbf{Ad}_{\mathbf{S}_{1,2}}{^{2}}\mathbf{V}$. The power conservation
yields that a wrench represented in $\mathcal{F}_{1}$ transforms to its
representation in $\mathcal{F}_{2}$ according to 
\begin{equation}
{^{2}}\mathbf{W}=\mathbf{Ad}_{\mathbf{S}_{1,2}}^{T}{^{1}}\mathbf{W}.
\end{equation}

While this notation is useful for kinetostatic formulations, it is
inconsistent in the sense that it treats screw coordinates differently for
twists and wrenches. In screw theory, aiming on a consistent treatment of
screw entities, a screw is represented by its coordinates as defined by (67)
in \cite{Part1} and the so-called \emph{reciprocal product} of two screws is
used \cite{Angeles2003,BottemaRoth1990,Selig}. The latter is defined for $%
\mathbf{X}_{1}=(\mathbf{\xi }_{1},\mathbf{\eta }_{1})^{T}$ and $\mathbf{X}%
_{2}=(\mathbf{\xi }_{2},\mathbf{\eta }_{2})^{T}$ as $\mathbf{X}_{1}\odot 
\mathbf{X}_{2}=\mathbf{\xi }_{1}^{T}\mathbf{\eta }_{2}+\mathbf{\eta }_{1}^{T}%
\mathbf{\xi }_{2}$. Two screws are said to be \emph{reciprocal} if $\mathbf{X%
}_{1}\odot \mathbf{X}_{2}=0$. Obviously, if twists and wrenches are
represented consistently with the same definition of screw coordinates, a
reciprocal twist and wrench screws means that they perform no work.
Geometrically, for zero pitch screws, this means that the screw axes
intersect.

In screw theory wrench screws are called co-screws to distinguish them from
motion screws and to indicate that a wrench acts on a motion screw (a twist)
as a linear operator that returns work or power. As twists form the Lie
algebra $se\left( 3\right) $, wrenches from the dual $se^{\ast }\left(
3\right) $.%
\newpage%

\section*{C Nomenclature}

\ 

\begin{tabular}{lllll}
$\mathcal{F}_{0}$ & - IFR &  &  &  \\ 
$\mathcal{F}_{i}$ & - BFR of body $i$ &  &  &  \\ 
$\mathcal{J}_{i,i}$ & - JFR for joint $i$ at body $i$, joint $i$ connects
body $i$ with its predecessor body $i-1$ &  &  &  \\ 
$\mathcal{J}_{i-1,i}$ & - JFR for joint $i$ at body $i-1$ &  &  &  \\ 
${^{i}}\mathbf{r}$ & - Coordinate representation of a vector resolved in BFR
on body $i$. &  &  &  \\ 
& \ \ The index is omitted if this is the IFR: $\mathbf{r}\equiv {^{0}}%
\mathbf{r}$. &  &  &  \\ 
$\mathbf{R}_{i}$ & - Rotation matrix from BFR $\mathcal{F}_{i}$ at body $i$
to IFR $\mathcal{F}_{0}$ &  &  &  \\ 
$\mathbf{R}_{i,j}$ & - Rotation matrix transforming coordinates resolved in
BFR $\mathcal{F}_{j}$ &  &  &  \\ 
& \ \ to coordinates resolved in $\mathcal{F}_{i}$ &  &  &  \\ 
$\mathbf{r}_{i}$ & - Position vector of origin of BFR $\mathcal{F}_{i}$ at
body $i$ resolved in IFR $\mathcal{F}_{0}$ &  &  &  \\ 
$\mathbf{r}_{i,j}$ & - Position vector from origin of BFR $\mathcal{F}_{i}$
to origin of BFR $\mathcal{F}_{j}$ &  &  &  \\ 
$\widetilde{\mathbf{x}}$ & - skew symmetric matrix associated to the vector $%
\mathbf{x}\in {\mathbb{R}}^{3}$ &  &  &  \\ 
$C_{i}=\left( \mathbf{R}_{i},\mathbf{r}_{i}\right) 
\hspace{-2ex}%
$ & - Absolute configuration of body $i$. This is denoted in matrix form
with $\mathbf{C}_{i}$ &  &  &  \\ 
$\mathbf{C}_{i,j}=\mathbf{C}_{i}^{-1}\mathbf{C}_{j}%
\hspace{-2ex}%
$ & - Relative configuration of body $j$ w.r.t. body $i$ &  &  &  \\ 
${^{k}}\mathbf{v}{_{i}^{j}}$ & - Translational velocity of body $i$ measured
at origin of BFR $\mathcal{F}_{j}$, resolved in BFR $\mathcal{F}_{k}$ &  & 
&  \\ 
$\mathbf{v}_{i}^{\text{b}}\equiv {^{i}}\mathbf{v}{_{i}^{i}}$ & - Body-fixed
representation of the translational velocity of body $i$ &  &  &  \\ 
$\mathbf{v}_{i}^{\text{s}}\equiv {^{0}}\mathbf{v}{_{i}^{0}}$ & - Spatial
representation of the translational velocity of body $i$ &  &  &  \\ 
${^{k}}\mathbf{\omega }{_{i}}$ & - Angular velocity of body $i$ measured and
resolved in BFR $\mathcal{F}_{k}$ &  &  &  \\ 
$\mathbf{\omega }{_{i}^{\text{b}}}\equiv {^{i}}\mathbf{\omega }{_{i}}$ & -
Body-fixed representation of the angular velocity of body $i$ &  &  &  \\ 
$\mathbf{\omega }{_{i}^{\text{s}}}\equiv {^{0}}\mathbf{\omega }{_{i}}$ & -
Spatial representation of the angular velocity of body $i$ &  &  &  \\ 
${^{k}}\mathbf{V}{_{i}^{j}}$ & - Twist of (RFR of) body $i$ measured in $%
\mathcal{F}_{j}$ and resolved in $\mathcal{F}_{k}$ &  &  &  \\ 
$\mathbf{V}_{i}^{\text{b}}\equiv {^{i}}\mathbf{V}{_{i}^{i}}$ & - Body-fixed
representation of the twist of body $i$ &  &  &  \\ 
$\mathbf{V}{_{i}^{\text{s}}}\equiv {^{0}}\mathbf{V}{_{i}^{0}}$ & - Spatial
representation of the twist of body $i$ &  &  &  \\ 
$\mathbf{V}_{i}^{\text{h}}\equiv {^{0}}\mathbf{V}{_{i}^{i}}$ & - Hybrid form
of the twist of body $i$ &  &  &  \\ 
$\mathsf{V}{^{\text{b}}}$ & - Vector of system twists in body-fixed
representation &  &  &  \\ 
$\mathsf{V}{^{\text{s}}}$ & - Vector of system twists in spatial
representation &  &  &  \\ 
$\mathsf{V}{^{\text{h}}}$ & - Vector of system twists in hybrid
representation &  &  &  \\ 
$\mathsf{V}{^{\text{m}}}$ & - Vector of system twists in mixed representation
&  &  &  \\ 
$\mathbf{W}_{i}^{\text{b}}$ & - Applied wrench at body $i$ in body-fixed
representation &  &  &  \\ 
$\mathbf{W}_{i}^{\text{s}}$ & - Applied wrench at body $i$ in spatial
representation &  &  &  \\ 
$\mathbf{W}_{i}^{\text{h}}$ & - Applied wrench at body $i$ in hybrid
representation &  &  &  \\ 
$\mathbf{M}_{i}^{\text{b}}$ & - Inertia matrix of body $i$ in body-fixed
representation &  &  &  \\ 
$\mathbf{M}_{i}^{\text{s}}$ & - Inertia matrix of body $i$ in spatial
representation &  &  &  \\ 
$\mathbf{M}_{i}^{\text{h}}$ & - Inertia matrix of body $i$ in hybrid
representation &  &  &  \\ 
$\mathbf{Ad}_{\mathbf{R}}$ & - Screw transformation associated with $%
C=\left( \mathbf{R},\mathbf{0}\right) $ &  &  &  \\ 
$\mathbf{Ad}_{\mathbf{r}}$ & - Screw transformation associated with $%
C=\left( \mathbf{I},\mathbf{r}\right) $ &  &  &  \\ 
$\mathbf{Ad}_{\mathbf{C}_{i,j}}$ & - Transformation matrix transforming
screw coordinates represented in $\mathcal{F}_{j}$ &  &  &  \\ 
& \ \ to screw coordinates represented in $\mathcal{F}_{i}$ &  &  &  \\ 
$\mathbf{ad}_{\mathbf{X}}$ & - Screw product matrix associated with screw
coordinate vector $\mathbf{X}\in {\mathbb{R}}^{6}$ &  &  &  \\ 
$\left[ \mathbf{X},\mathbf{Y}\right] $ & - Lie bracket of screw coordinate
vectors $\mathbf{X},\mathbf{Y}\in {\mathbb{R}}^{6}$. It holds $\left[ 
\mathbf{X},\mathbf{Y}\right] =\mathbf{ad}_{\mathbf{X}}\mathbf{Y}$. &  &  & 
\\ 
$\widehat{\mathbf{X}}\in se\left( 3\right) $ & - $4\times 4$ matrix
associated with the screw coordinate vectors $\mathbf{X}\in {\mathbb{R}}^{6}$
&  &  &  \\ 
$SE\left( 3\right) $ & - Special Euclidean group in three dimensions -- Lie
group of rigid body motions &  &  &  \\ 
$se\left( 3\right) $ & - Lie algebra of $SE\left( 3\right) $ -- algebra of
screws &  &  &  \\ 
$\mathbf{q}\in {\mathbb{V}}^{n}$ & - Joint coordinate vector &  &  &  \\ 
${\mathbb{V}}^{n}$ & - Configuration space &  &  & 
\end{tabular}

\section*{Acknowledgement}

The author acknowledges that this work has been partially supported by the
Austrian COMET-K2 program of the Linz Center of Mechatronics (LCM).

\end{document}